\documentclass[11pt]{article}

\usepackage{amsmath}
\usepackage{amssymb}
\usepackage{amsfonts}
\usepackage{amsthm}
\usepackage{array}

\usepackage[small,nohug,heads=vee]{diagrams}
\diagramstyle[labelstyle=\scriptstyle]
\usepackage{diagrams}

\newcolumntype{C}[1]{>{\centering\hspace{0pt}}p{#1}}
\usepackage[pdftex]{graphicx}
\usepackage[margin=1in]{geometry}

\usepackage[nottoc]{tocbibind}

\title{Nearly-K\"{a}hler $6$-Manifolds of Cohomogeneity Two: Local Theory}
\author{Jesse Madnick}
\date{October 2018}

\newcommand{\Addresses}
{{  \bigskip
  \textsc{McMaster University} \par\nopagebreak
  \textsc{Department of Mathematics \& Statistics}\par\nopagebreak
  \textsc{Hamilton, ON, Canada, L8S 4K1}\par\nopagebreak
  \textit{E-mail address}: \texttt{madnickj@mcmaster.ca}
}}

\begin{document}

\maketitle

\begin{abstract}
We study nearly-K\"{a}hler $6$-manifolds equipped with a cohomogeneity-two Lie group action for which the principal orbits are coisotropic.  If the metric is complete, then we show that this last condition is automatically satisfied, and both the acting Lie group and the principal orbits are finite quotients of $\mathbb{S}^3 \times \mathbb{S}^1$.  \\
\indent We partition the class of such nearly-K\"{a}hler structures into three types (called I, II, III) and prove a local existence and generality result for each type.  Metrics of Types I and II are shown to be incomplete. \\
\indent We also derive a quasilinear elliptic PDE system on the 2-dimensional orbit space which nearly-K\"{a}hler structures of Type I must satisfy.  Finally, we remark on a relatively simple one-parameter family of Type III structures that turn out to be incomplete metrics that are cohomogeneity-one under the action of a larger group.
\end{abstract}

\tableofcontents

\pagebreak
 
\section{Introduction}
   
\indent \indent Nearly-K\"{a}hler $6$-manifolds are a class of Riemannian $6$-manifolds $(M^6, g)$ whose geometry is in some sense modeled on the round $6$-sphere $\mathbb{S}^6 \simeq \text{G}_2/\text{SU}(3)$.  Like the round $\mathbb{S}^6$, they carry a triple $(J, \Omega, \Upsilon)$ consisting of a compatible almost-complex structure $J$, a non-degenerate $2$-form $\Omega$, and a $(3,0)$-form $\Upsilon$, and these are asked to satisfy the defining differential equations
\begin{align*}
d\Omega & = 3\,\text{Im}(\Upsilon) \\
d\,\text{Re}(\Upsilon) & = 2\,\Omega \wedge \Omega.
\end{align*}
Here, the almost-complex structure $J$ is \textit{not} integrable, and the $2$-form $\Omega$ is \textit{not} closed. \\
\indent Yet, in spite of these two shortcomings, nearly-K\"{a}hler $6$-manifolds enjoy several remarkable characterizations that have led to increased attention as of late, especially in connection with exceptional holonomy metrics and real Killing spinors.  Indeed: \\

\noindent \textbf{Theorem (B\"{a}r \cite{MR1224089}, Grunewald \cite{Grunewald1990}):} Let $(M^6, g)$ be a simply-connected, spin, complete Riemannian $6$-manifold.  Let $\text{Cone}(M,g) = (\mathbb{R}^+ \times M, dt^2 + t^2g)$ be its Riemannian cone.  The following are equivalent: \\
\indent (i) $(M, g)$ admits a nearly-K\"{a}hler structure. \\
\indent (ii) $(M, g)$ has a real Killing spinor. \\
\indent (iii) $\text{Cone}(M,g)$ has a parallel spinor. \\
\indent (iv) $\text{Hol}(\text{Cone}(M,g)) = \text{G}_2$ or $(M^6, g) \cong (\mathbb{S}^6, g_{\text{round}})$. \\

\indent In fact, nearly-K\"{a}hler $6$-manifolds are Einstein of positive scalar curvature.  Thus, complete examples are compact with finite fundamental group (by Bonnet-Myers). \\

\indent A central problem in the study of nearly-K\"{a}hler $6$-manifolds is the present dearth of compact, simply-connected examples.  Indeed, as of this writing, only six such examples are known.  Four of these are the homogeneous spaces \cite{MR0236329}
$$\mathbb{S}^6 = \frac{\text{G}_2}{\text{SU}(3)}, \ \ \ \ \ \mathbb{S}^3 \times \mathbb{S}^3 = \frac{\text{SU}(2)^3}{\Delta\text{SU}(2)}, \ \ \ \ \ \text{CP}^3 = \frac{\text{Sp}(2)}{\text{U}(1) \times \text{Sp}(1)}, \ \ \ \ \ \text{Flag}(\mathbb{C}^3) = \frac{\text{SU(3)}}{T^2},$$
and it has been shown \cite{MR2158165} that these are the only possible homogeneous examples.  Here, we caution that the metric on $\mathbb{S}^3 \times \mathbb{S}^3$ is \textit{not} the product metric, and the almost-complex structure on $\text{CP}^3$ is \textit{not} the standard one. \\
\indent Following work of Conti and Salamon \cite{MR2327032}, Fern\'{a}ndez, Ivanov, Mu\~noz and Ugarte \cite{MR2456893}, and Podest\`{a} and Spiro \cite{MR2587385} \cite{MR2917173}, recently Foscolo and Haskins \cite{MR3583352} succeeded in constructing \textit{inhomogeneous} nearly-K\"{a}hler metrics on $\mathbb{S}^6$ and $\mathbb{S}^3 \times \mathbb{S}^3$ that are cohomogeneity-one under an $(\text{SU}(2) \times \text{SU}(2))$-action.  Their approach involves cohomogeneity-one techniques, drawing on methods of Eschenburg and Wang \cite{MR1758585} and B\"{o}hm \cite{MR1700472}, guided by the idea that such examples might arise as desingularizations of the sine-cone over the Sasaki-Einstein $\mathbb{S}^2 \times \mathbb{S}^3$. \\
\indent From the point of view of symmetries, the next natural question is the existence of compact simply-connected examples of cohomogeneity-two.  This remains a difficult open problem, and is the primary motivation for this work.

\subsection{Methods and Main Results}

\indent \indent In this work, we study the geometry of cohomogeneity-two nearly-K\"{a}hler $6$-manifolds.  That is, we study nearly-K\"{a}hler $6$-manifolds $M$ equipped with a faithful $G$-action whose generic orbits have codimension two.  We always suppose that $G$ is a closed, connected subgroup of the isometry group of $M$, and that the $G$-action preserves $(J, \Omega, \Upsilon)$. \\
\indent We will restrict attention to the case where the principal $G$-orbits are \textit{coisotropic}, meaning that the $4$-form $\Omega \wedge \Omega$ vanishes on these (4-dimensional) orbits.  Our first result shows that, in fact, this case is the one of most interest: \\

\noindent \textbf{Theorem 1.1:} Let $M$ be a nearly-K\"{a}hler $6$-manifold.  Suppose that a connected Lie group $G$, closed in the isometry group of $M$, acts faithfully on $M$ with cohomogeneity two, preserving $(J, \Omega, \Upsilon)$. \\
\indent (a) If $M$ is complete, then the principal $G$-orbits in $M$ are coisotropic. \\
\indent (b) If the principal $G$-orbits are coisotropic, then $G$ is $4$-dimensional and non-abelian. \\
\indent (c) If $M$ is complete, then both $G$ and the principal $G$-orbits in $M$ are finite quotients of $\mathbb{S}^3 \times \mathbb{S}^1$. \\

\indent Next, we turn to the question of local existence.  That is, on sufficiently small open sets of $\mathbb{R}^6$ we ask whether cohomogeneity-two nearly-K\"{a}hler metrics can exist at all.  If so, what is the initial data required to construct these metrics as solutions to a (sequence of) Cauchy problem(s)? \\
\indent We approach this problem by an application of Cartan's Third Theorem \cite{Bryant:2014aa}.  This result generalizes Lie's Third Theorem on the ``integration" of Lie algebras to local Lie groups.  Its primary hypothesis is that ``mixed partials commute," meaning the satisfaction of a set of integrability conditions (analogous to the Jacobi identity for Lie algebras). \\
\indent In the case of cohomogeneity-two nearly-K\"{a}hler metrics with coisotropic principal orbits, these integrability conditions form a system of 80 quadratic equations for 55 unknown functions.  Careful study of this system leads us to partition the class of metrics under consideration into three types, called Types I, II, and III. \\
\indent We will show that metrics of each Type exist locally and in abundance: each Type is an infinite-dimensional family.  More precisely: \\

\noindent \textbf{Theorem 1.2:} On sufficiently small open sets in $\mathbb{R}^6$: \\
\indent (a) Nearly-K\"{a}hler structures of Type I exist, depending on 2 arbitrary functions of $1$ variable.  If $M$ is of Type I, then $G$ is a discrete quotient of $\text{H}_3 \times \mathbb{R}$, where $\text{H}_3$ is the real Heisenberg group.  In particular, metrics of Type I are incomplete. \\
\indent (b) Nearly-K\"{a}hler structures of Type II exist, depending on 2 arbitrary functions of $1$ variable.  If $M$ is of Type II, then $G$ is solvable.  In particular, metrics of Type II are incomplete. \\
\indent (c) Nearly-K\"{a}hler structures of Type III with $G = (\mathbb{S}^3 \times \mathbb{S}^1)/(\text{Finite})$ exist, depending on $2$ arbitrary functions of $1$ variable. \\

\indent The dependence on $2$ arbitrary functions of $1$ variable --- the same initial data (or ``local generality") required to construct holomorphic functions $f \colon \mathbb{C} \to \mathbb{C}$ --- suggests the possibility that cohomogeneity-two nearly-K\"{a}hler $6$-manifolds may be recovered from holomorphic data.  More precisely, one can ask:
\begin{enumerate}
\item Can cohomogeneity-two nearly-K\"{a}hler structures be reconstructed from solutions to an elliptic PDE system on a Riemann surface $\Sigma$?
\item Can solutions to this elliptic PDE system, in turn, be reinterpreted as pseudo-holomorphic curves in some almost-complex manifold?
\item Can this elliptic PDE system be recast as a single second-order elliptic PDE on $\Sigma$?
\item Can cohomogeneity-two nearly-K\"{a}hler structures be reconstructed from holomorphic data by means of a Weierstrass representation formula?
\end{enumerate}

\indent These questions provide interesting directions for further study.  In this work, we offer the following first steps in the Type I setting: \\

\noindent \textbf{Proposition 1.3:} Across their principal loci, nearly-K\"{a}hler structures of Type I are solutions to a certain quasilinear elliptic PDE system (5.11) on a Riemann surface. \\

\indent The derivation of the PDE system (5.11) can be mimicked to obtain similar elliptic PDE systems in the Type II and Type III settings as well.  However, as we hope to make plain, the PDE systems in those settings are extremely cumbersome to write down explicitly. \\

\indent We expect an affirmative answer to Question 2, and are optimistic about an affirmative answer to Question 3, as well.  Of course, Question 4 is significantly more speculative. \\

\indent At the end of this report, we study a particular subclass of the Type III nearly-K\"{a}hler metrics with $G = (\mathbb{S}^3 \times \mathbb{S}^1)/(\text{Finite})$ in some detail.  This subclass turns out to be a one-parameter family of metrics $\{g_c \colon c \in [\sqrt{3}, \infty]\}$, and we derive an explicit formula for $g_\infty$.  Unfortunately, the metric $g_{\sqrt{3}}$ is homogeneous, while the metrics $g_c$ for $c \in (\sqrt{3}, \infty]$ are incomplete and cohomogeneity-one under the action of a larger group.

\subsection{Organization}

\indent \indent This work is organized as follows.  In $\S2$, we review the fundamentals of $H$-structures and intrinsic torsion, the language in which this work is phrased.  In particular, we state Cartan's Third Theorem (labeled Theorem 2.2), our main tool for proving local existence/generality results. \\
\indent In $\S3.1$, we compare and contrast various definitions of ``nearly-K\"{a}hler $6$-manifold" encountered in the literature, clarifying our own conventions.  The material of $\S2$ and $\S3.1$ is  standard, and experts may wish to skip these. In $\S3.2$, we prove Theorem 1.1(a) (labeled Proposition 3.3).  \\

\indent Section 4 sets up the moving frame apparatus that we will use in our study.  In $\S$4.2, we adapt frames to the (coisotropic) principal $G$-orbits.  This frame adaptation defines an $\text{O}(2)$-structure, and its study is central to this work.  In $\S$4.3, we describe the intrinsic torsion of our $\text{O}(2)$-structure,  concrete geometric interpretations of which are offered in $\S$4.5. \\
\indent In $\S$4.4, we prove Theorem 1.1(b) and 1.1(c) (labeled Proposition 4.5).  Section 4.6 sets up the machinery needed to describe cohomogeneity-two nearly-K\"{a}hler structures as solutions to an elliptic PDE system on a Riemann surface. \\

\indent In $\S$5, we describe our partition into Types I, II, and III.  In $\S$5.1, we examine Type I structures and prove Theorem 1.2(a) (labeled Theorem 5.4 and Proposition 5.5) and Proposition 1.3 (labeled Proposition 5.6).  Similarly, $\S$5.2 pertains to Type II structures and contains a proof of Theorem 1.2(b) (labeled Theorem 5.9 and Proposition 5.10), and $\S$5.3 contains a proof of Theorem 1.2(c) (labeled Theorem 5.15). \\
\indent Finally, in $\S$6, we study the one-parameter family of metrics $\{g_c \colon c \in [\sqrt{3}, \infty]\}$ mentioned above. \\

\noindent \textbf{Notation:} The following notation and terminology will be used throughout.
\begin{itemize}
  \item Let $\pi \colon P \to M$ be a submersion.  A $k$-form $\theta \in \Omega^k(P)$ is \textit{$\pi$-semibasic} if $X\,\lrcorner\,\theta = 0$ for all vectors $X \in TP$ tangent to the $\pi$-fibers.  We will simply say ``semibasic" when $\pi$ is clear from context.
  
  \item When $\omega = (\omega^1, \ldots, \omega^n)$ denotes the tautological $1$-form on an $H$-structure $B \to M^n$, we will use the shorthand
$$\omega^{ij} := \omega^i \wedge \omega^j, \ \ \ \ \ \omega^{ijk} = \omega^i \wedge \omega^j \wedge \omega^k, \ \ \ \ \ \text{etc.}$$
to denote wedge products.

  \item For $1$-forms $\alpha_1, \ldots, \alpha_k \in \Omega^1(M)$, we let $\langle \alpha_1, \ldots, \alpha_k \rangle$ denote the differential ideal in $\Omega^*(M)$ generated by these $1$-forms.  In particular, $\langle \alpha_1, \alpha_2 \rangle$ denotes an ideal (not an inner product). \\
\end{itemize}

\noindent \textbf{Acknowledgements:} This work forms part of the author's 2018 Ph.D. thesis at Stanford University. I would like to thank Robert Bryant for suggesting this problem and generously sharing his many insights.  Without his tireless guidance, this work would not have been possible.  I would also like to thank Rick Schoen for supporting my research throughout the duration of this project, and for all the inspiring lectures, seminars, and discussions on geometric analysis. \\
\indent This work has also benefited from helpful conversations with Gavin Ball, Jason DeVito, Christos Mantoulidis, Rafe Mazzeo, Gon\c calo Oliveira, McKenzie Wang, and Wolfgang Ziller.
    
\section{$H$-Structures and Cartan's Third Theorem}

\indent \indent Much of this work will be phrased in the language of $H$-structures, intrinsic torsion, and augmented coframings.  As such, we use this section to recall this terminology, set notation, and describe our primary technical tool for proving local existence.  The material in this section is standard; more information can be found in \cite{Bryant:2014aa}, \cite{MR1062197}, and \cite{MR1004008}.

\subsection{$H$-Structures and Intrinsic Torsion}

\indent \indent Let $M$ be a smooth $n$-manifold.  A \textit{coframe} at $x \in M$ is a vector space isomorphism $u \colon T_xM \to \mathbb{R}^n$.  We let $\pi \colon FM \to M$ denote the \textit{general coframe bundle}, which is the principal right $\text{GL}_n(\mathbb{R})$-bundle over $M$ whose fiber at $x \in M$ consists of the coframes at $x$.  Here, the right $\text{GL}_n(\mathbb{R})$-action on $FM$ is by composition: for $g \in \text{GL}_n(\mathbb{R})$ and $u \in FM$, we set
$$u \cdot g := g^{-1} \circ u.$$
\indent A \textit{coframing} on an open set $U \subset M$ is an $n$-tuple $\eta = (\eta^1, \ldots, \eta^n)$ of linearly independent $1$-forms on $U$.  We think of coframings as $\mathbb{R}^n$-valued $1$-forms $\eta \in \Omega^1(U; \mathbb{R}^n)$ for which each $\eta_x \colon T_xU \to \mathbb{R}^n$ is a coframe.  Alternatively, we regard coframings as local sections $\sigma_\eta \in \Gamma(U; FM)$ or as local trivializations $\psi_\eta \colon U \times \text{GL}_n(\mathbb{R}) \to FM|_U$ via $\psi_\eta(x,g) = \eta_x \cdot g$. \\
\indent To a local diffeomorphism $f \colon M_1 \to M_2$, we associate the bundle map $f^{(1)} \colon FM_1 \to FM_2$ defined by
$$f^{(1)}(u) = u \circ (f_*|_{\pi_1(u)})^{-1}.$$
One can check that $f \mapsto f^{(1)}$ is functorial. \\

\indent For a subgroup $H \leq \text{GL}_n(\mathbb{R})$, an \textit{$H$-structure} $B$ on an $n$-manifold $M^n$ is an $H$-subbundle of the general coframe bundle $B \subset FM$.  Note that, despite the terminology, an $H$-structure depends on the representation of $H$ on $\mathbb{R}^n$, not just on the abstract group itself. \\
\indent We say that $H$-structures $\pi_1 \colon B_1 \to M$ and $\pi_2 \colon B_2 \to M_2$ are \textit{(locally) equivalent} if there is a (local) diffeomorphism $f \colon M_1 \to M_2$ for which $f^{(1)}(B_1) = B_2$. \\

\indent The \textit{tautological $1$-form} on an $H$-structure $B$ is the $\mathbb{R}^n$-valued $1$-form $\omega = (\omega^1, \ldots,$ $\omega^n) \in \Omega^1(B; \mathbb{R}^n)$ given by
$$\omega(v) = u(\pi_*(v)), \ \ \text{for } v \in T_uB.$$
The tautological $1$-form ``reproduces" all of the local coframings of $M$, in that it satisfies the following property: For any coframing $\eta \in \Omega^1(U; \mathbb{R}^n)$, we have $\sigma_\eta^*(\omega^1, \ldots, \omega^n) = (\eta^1, \ldots, \eta^n)$, or equivalently, $\psi_\eta^*(\omega^1, \ldots, \omega^n)|_{(x, h)} = (\eta^1, \ldots, \eta^n)|_x \cdot h$. \\
\indent One can show \cite{MR1062197} that if $H$ is connected, a smooth map $F \colon B_1 \to B_2$ between $H$-structures is a local equivalence if and only if $F^*(\omega_2) = \omega_1$. \\

\indent A \textit{connection} on an $H$-structure $B$ is simply a connection on the principal $H$-bundle $B$.  That is, it is an $\mathfrak{h}$-valued $1$-form $\phi \in \Omega^1(B; \mathfrak{h})$ that sends $H$-action vector fields to their Lie algebra generators and is $H$-equivariant:
\begin{align*}
\phi(X^\#) & = X, \ \text{ for all } X \in \mathfrak{h} \\
R_h^*(\phi) & = \text{Ad}_{h^{-1}}(\phi), \ \text{ for all } h \in H.
\end{align*}
Note that the first condition implies that $\phi$ restricts to each $H$-fiber to be the Maurer-Cartan form on $H$. \\
\indent Given an $H$-structure $\pi \colon B \to M$ with connection $\phi \in \Omega^1(B; \mathfrak{h})$, one can differentiate the equation $\psi_\eta^*(\omega) = \eta \cdot h$ to derive \textit{Cartan's first structure equation}
$$d\omega = -\phi \wedge \omega + \textstyle \frac{1}{2}T_\phi(\omega \wedge \omega),$$
where $T_\phi \colon B \to \mathbb{R}^n \otimes \Lambda^2(\mathbb{R}^n)^*$ is a function called the \textit{torsion of the connection} $\phi$.  To emphasize the distinction between $\mathbb{R}^n$ and $(\mathbb{R}^n)^*$, let us write $V = \mathbb{R}^n$. \\

\indent Let $\phi_1$, $\phi_2$ be two connections on $B$, with torsion functions $T_{\phi_1}$, $T_{\phi_2}$, respectively.  The difference $\phi_1 - \phi_2$ is $\pi$-semibasic and so can be written $\phi_1 - \phi_2 = p(\omega)$ for some function $p \colon B \to \mathfrak{h} \otimes V^*$.  A calculation shows \cite{MR1062197}, \cite{MR1004008} that the difference in the torsion functions is
$$T_{\phi_1} - T_{\phi_2} = \delta(p),$$
where $\delta \colon \mathfrak{h} \otimes V^* \hookrightarrow V \otimes V^* \otimes V^* \to V \otimes \Lambda^2V^*$ is the $H$-equivariant linear map given by skew-symmetrization.  Thus, the composite map
$$T \colon B \to V \otimes \Lambda^2V^* \twoheadrightarrow \frac{V \otimes \Lambda^2V^*}{\delta(\mathfrak{h} \otimes V^*)} =: H^{0,2}(\mathfrak{h})$$
is independent of the choice of connection $\phi$.  We refer to $T$ as the \textit{intrinsic torsion of the $H$-structure}, and the codomain $H^{0,2}(\mathfrak{h}) = (V \otimes \Lambda^2V^*)/\text{Im}(\delta)$ as the \textit{intrinsic torsion space}. \\

\noindent \textit{Remark:} The vector space $H^{0,2}(\mathfrak{h})$ can be regarded as a Spencer cohomology group, which explains the reason for the notation. $\Box$

\subsection{The Case of $H \leq \text{SO}(n)$}

\indent  \indent Suppose now that $H \leq \text{SO}(n)$.  We regard $B \subset F_{\text{SO}(n)}$, where $F_{\text{SO}(n)}$ is the orthonormal frame bundle corresponding to the underlying $\text{SO}(n)$-structure.  Let $\theta \in \Omega^1(F_{\text{SO}(n)}; \mathfrak{so}(n))$ denote the Levi-Civita connection.  On $F_{\text{SO}(n)}$, the Fundamental Lemma of Riemannian Geometry gives
$$d\omega = -\theta \wedge \omega.$$
Let us split $\mathfrak{so}(n) = \mathfrak{h} \oplus \mathfrak{h}^\perp$ with respect to the Killing form of $\mathfrak{so}(n)$.  Accordingly, we split
\begin{align*}
\theta|_B = \gamma_H + \tau_H,
\tag{2.1}
\end{align*}
where $\gamma_H \in \Omega^1(B; \mathfrak{h})$ and $\tau_H \in \Omega^1(B; \mathfrak{h}^\perp)$.  One can check that $\gamma_H$ is a connection on the $H$-structure $B$, while $\tau_H = t(\omega)$ for some $t \colon B \to \mathfrak{h}^\perp \otimes V^*$.  Thus, on $B$,
\begin{align*}
d\omega & = \textstyle -\gamma_H \wedge \omega + \frac{1}{2}\delta(t)(\omega \wedge \omega),
\end{align*}
and so the torsion of the connection $\gamma_H$ takes values in $\delta(\mathfrak{h}^\perp \otimes V^*)$.  In fact, since $\delta \colon \mathfrak{so}(n) \otimes V^* \to V \otimes \Lambda^2V^*$ is injective, and since $\Lambda^2V^* \cong \mathfrak{so}(n) = \mathfrak{h} \oplus \mathfrak{h}^\perp$, we have
\begin{align*}
V \otimes \Lambda^2V^* \cong  V \otimes (\mathfrak{h} \oplus \mathfrak{h}^\perp) & = (V \otimes \mathfrak{h}) \oplus (V \otimes \mathfrak{h}^\perp) \cong \delta(\mathfrak{h} \otimes V) \oplus \delta(\mathfrak{h}^\perp \otimes V),
\end{align*}
whence
$$H^{0,2}(\mathfrak{h}) \cong \delta(\mathfrak{h}^\perp \otimes V) \cong \mathfrak{h}^\perp \otimes V.$$
We will return to this formula in $\S$3.1 in the cases $H = \text{U}(3) \leq \text{SO}(6)$ and $H = \text{SU}(3)  \leq \text{SO}(6)$.

\subsection{Group Actions on $H$-Structures}

\indent \indent We will be concerned with $H$-structures on manifolds $M$ equipped with a $G$-action that preserves the $H$-structure.  In this regard, we make a simple preliminary observation. \\
\indent A $G$-action on $M$ induces $G$-actions on both $T^*M$ and $FM$.  Explicitly, the $G$-action on $FM$ is
$$g \cdot u = (g^{-1})^*u = u \circ (g^{-1})_*.$$
Note that if $g \in G$ stabilizes a coframe $u \in FM|_x$, then $gx = x$ and $(g^{-1})^*u = u$, so that $g$ acts as the identity on $T^*_xM$.  From this, we observe: \\

\noindent \textbf{Lemma 2.1:} Let $P \to M^n$ be an $H$-structure, where $H \leq \text{SO}(n)$.  Suppose $M$ is equipped with a $G$-action that preserves the $H$-structure and acts by cohomogeneity-$k$ on $M$.  Then $n - k \leq \dim(G) \leq n + \dim(H)$. \\

\noindent \textit{Proof:} Since $G$ acts with cohomogeneity-$k$ on $M^n$, so $G$ acts transitively on the $(n-k)$-dimensional principal orbits in $M$, so $\dim(G) \geq n-k$. \\
\indent On the other hand, if $g$ stabilizes a coframe $u \in FM|_x$, then $g$ acts as the identity on $T^*_xM$.  Since $g$ is an isometry (because $H \leq \text{SO}(n)$), so $g = \text{Id}$, so the $G$-action on $P$ is free.  Thus, $\dim(G) \leq \dim(P) = n + \dim(H)$. $\lozenge$

\subsection{Cartan's Third Theorem}

\indent \indent In order to prove the local existence of $H$-structures with desired properties, we encode the data of an $H$-structure in terms of an ``augmented coframing." \\

\noindent \textbf{Definition:} An \textit{augmented coframing} on an $n$-manifold $P$ is a triple $(\eta, a,b)$, where $\eta = (\eta^1, \ldots, \eta^n)$ is a coframing on $P$, and $a = (a^1, \ldots, a^s) \colon P \to \mathbb{R}^s$ and $b = (b^1, \ldots, b^r) \colon P \to \mathbb{R}^r$ are smooth functions. \\
\indent The functions $a^1, \ldots, a^s \colon P \to \mathbb{R}$ are called the \textit{primary invariants} of the augmented coframing, while the functions $b^1, \ldots, b^r \colon P \to \mathbb{R}$ are called \textit{free derivatives}. \\
\indent For the rest of this section, we fix index ranges $1 \leq i,j,k \leq n$ and $1 \leq \alpha, \beta \leq s$ and $1 \leq \rho \leq r$.

\pagebreak

\indent We will be interested in augmented coframings that satisfy a given set of \textit{structure equations}, by which we mean a set of equations of the form
\begin{align*}
d\eta^i & = \textstyle -\frac{1}{2}C^i_{jk}(a)\,\eta^j \wedge \eta^k \tag{2.2} \\
da^\alpha & = F^\alpha_i(a,b)\,\eta^i
\end{align*}
for some given functions $C^i_{jk}(u) = -C^i_{jk}(u)$ on $\mathbb{R}^s$ and $F^\alpha_i(u,v)$ on $\mathbb{R}^s \times \mathbb{R}^r$. \\

\indent Let $\pi \colon B \to M^n$ and $\theta \in \Omega^1(B; \mathfrak{h})$ be an $H$-structure-with-connection.  Let $\omega \in \Omega^1(B; \mathbb{R}^n)$ denote the tautological $1$-form on $B$.  Then $\eta = (\omega, \theta) = (\omega^i, \theta^j_k) \colon TB \to \mathbb{R}^n \oplus \mathfrak{h}$ is a coframing of $B$ whose exterior derivatives satisfy equations of the form
\begin{align*}
d\omega^i & = -\theta^i_j \wedge \omega^j + T^i_{jk}\,\omega^j \wedge \omega^k \tag{2.3a} \\
d\theta^i_j & = -\theta^i_k \wedge \theta^k_j + R^i_{jk\ell}\,\omega^k \wedge \omega^\ell \tag{2.3b} \\
dT^i_{jk} & = A^i_{jk\ell}(T)\, \theta^\ell + B^i_{jk\ell} \omega^\ell \tag{2.3c} \\
dR^i_{jk\ell} & = C^i_{jk\ell m}(R)\, \theta^m + D^i_{jk\ell m} \omega^m. \tag{2.3d}
\end{align*}
for some functions $T = (T^i_{jk}) \colon B \to V \otimes \Lambda^2V^*$ and $R = (R^i_{jk\ell}) \colon B \to \mathfrak{h} \otimes \Lambda^2V^*$. \\
\indent Conversely, suppose $P$ is a manifold with a coframing $\eta = (\omega, \theta) \colon TP \to \mathbb{R}^n \oplus \mathfrak{h}$ and functions $T = (T^i_{jk}) \colon P \to V \otimes \Lambda^2V^*$ and $R = (R^i_{jk\ell}) \colon P \to \mathfrak{h} \otimes \Lambda^2V^*$ satisfying (2.3a)-(2.3d).  From (2.3a), there is a submersion $\pi \colon P \to M$ whose fibers are integral manifolds of the (Frobenius) ideal $\langle \omega^1, \ldots, \omega^n \rangle$.  Further, one can construct a local diffeomorphism $\sigma \colon P \to FM$ whose image is an $H$-structure $B \subset FM$ such that $\sigma$ sends $\pi$-fibers to $H$-orbits and has $\sigma^*(\omega_0) = \omega$, where $\omega_0$ is the tautological form on $B$. \\

\indent To prove the local existence of augmented coframings satisfying prescribed structure equations (2.2), we will appeal to a very general result.  This theorem, due to Cartan, is a vast generalization of the converse to Lie's Third Theorem on the ``integration" of a Lie algebra to a local Lie group.  Roughly, the theorem says that the necessary first-order conditions for existence --- namely, $d(d\eta^i) = 0$ and $d(da^\alpha) = 0$ --- are very close to sufficient. \\
\indent Let us be more explicit.  The equations $d(d\eta^i) = 0$, meaning $d(C^i_{jk}(a)\,\eta^j \wedge \eta^k) = 0$, expand to
\begin{equation*}
F^\alpha_j \frac{\partial C^i_{k\ell}}{\partial u^\alpha} + F^\alpha_k \frac{\partial C^i_{\ell j}}{\partial u^\alpha} + F^\alpha_\ell \frac{\partial C^i_{jk}}{\partial u^\alpha} = C^i_{mj}C^m_{k\ell} + C^i_{mk}C^m_{\ell j} + C^i_{m\ell} C^m_{jk}.
\tag{2.4}
\end{equation*}
Similarly, the equations $d(da^\alpha) = 0$, meaning $d(F^\alpha_i(a,b)\,\eta^i) = 0$, expand to
\begin{align*}
0 & = \frac{\partial F^\alpha_i}{\partial v^\rho}\,db^\rho \wedge \eta^i + \frac{1}{2}\!\left( F^\beta_i \frac{\partial F^\alpha_j}{\partial u^\beta} - F^\beta_j \frac{\partial F^\alpha_i}{\partial u^\beta} - C^\ell_{ij}F^\alpha_\ell \right) \eta^i \wedge \eta^j.
\end{align*}
Since we lack formulas for $db^\rho$, it is not immediately clear how to satisfy this condition.  However, if there exist functions $G^\rho_j$ on $\mathbb{R}^s \times \mathbb{R}^r$ for which
\begin{equation*}
F^\beta_i \frac{\partial F^\alpha_j}{\partial u^\beta} - F^\beta_j \frac{\partial F^\alpha_i}{\partial u^\beta} - C^\ell_{ij} F^\alpha_\ell = \frac{\partial F^\alpha_i}{\partial v^\rho} G^\rho_j - \frac{\partial F^\alpha_j}{\partial v^\rho} G^\rho_i,
\tag{2.5}
\end{equation*}
then $d(da^\alpha) = 0$ reads simply
\begin{equation*}
0 = \frac{\partial F^\alpha_i}{\partial v^\rho} \left(db^\rho - G^\rho_j \,\eta^j\right) \wedge \eta^i.
\end{equation*}
Thus, if functions $G^\rho_j$ exist which satisfy (2.5), then there will exist an expression of the $db^\rho$ in terms of $\eta^i$ that will fulfill $d(da^\alpha) = 0$. We need one last piece of terminology before stating the theorem. \\

\noindent \textbf{Definition:} The \textit{tableau of free derivatives} of the equations (2.2) at a point $(u,v) \in \mathbb{R}^s \times \mathbb{R}^r$ is the linear subspace $A(u,v) \subset \text{Hom}(\mathbb{R}^n, \mathbb{R}^s)$ given by
$$A(u,v) = \text{span}\left\{ \frac{\partial F^\alpha_i}{\partial v^\rho}(u,v)\,e_\alpha \otimes f^i \colon 1 \leq \rho \leq r \right\},$$
where here $\{e_1, \ldots, e_s\}$ is a basis of $\mathbb{R}^s$ and $\{f^1, \ldots, f^n\}$ is a basis of $(\mathbb{R}^n)^*$. \\

\noindent \textbf{Theorem 2.2 (Cartan):} Fix real-analytic functions $C^i_{jk} = -C^k_{jk}$ on $\mathbb{R}^s$ and $F^\alpha_i$ on $\mathbb{R}^s \times \mathbb{R}^r$.  Suppose that: \\
\indent $\bullet$ The functions $C^i_{jk}$ and $F^\alpha_i$ satisfy (2.4). \\
\indent $\bullet$ There exist real-analytic functions $G^\rho_i$ on $\mathbb{R}^s \times \mathbb{R}^r$ satisfying (2.5). \\
\indent $\bullet$ The tableau of free derivatives $A(u,v)$ is involutive, has dimension $r$, and has Cartan characters $(\widetilde{s}_1, \ldots, \widetilde{s}_n)$ for all $(u,v) \in \mathbb{R}^s \times \mathbb{R}^r$. \\
\indent Then for any $(a_0, b_0) \in \mathbb{R}^s \times \mathbb{R}^r$, there exists a real-analytic augmented coframing $(\eta, a,b)$ on an open neighborhood of $0 \in \mathbb{R}^n$ that satisfies (2.2) and has $(a(0), b(0)) = (a_0, b_0)$. \\
\indent Moreover, augmented coframings satisfying (2.2) depend (modulo diffeomorphism) on $\widetilde{s}_p$ functions of $p$ variables (in the sense of exterior differential systems) where $\widetilde{s}_p$ is the last non-zero Cartan character of $A(u,v)$. \\
\\
\noindent \textit{Remark:} In outline, the proof of Theorem 2.2 is as follows: One constructs an exterior differential system on the manifold $\text{GL}_n(\mathbb{R}) \times \mathbb{R}^n \times \mathbb{R}^s \times \mathbb{R}^r$ whose integral $n$-manifolds are in bijection with augmented coframings satisfying (2.2).  An application of the Cartan-K\"{a}hler Theorem then yields the desired integral $n$-manifolds, and these depend on $\widetilde{s}_p$ functions of $p$ variables.  For details, see \cite{Bryant:2014aa}. \\
\indent The Cartan-K\"{a}hler Theorem requires real-analyticity, which is why Theorem 2.2 does, too.  However, since we will be using Theorem 2.2 to construct Einstein metrics --- which are real-analytic in harmonic coordinates \cite{MR644518} --- the real-analyticity hypothesis is not a significant limitation. $\Box$

\section{Nearly-K\"{a}hler $6$-Manifolds}

\subsection{Nearly-K\"{a}hler $6$-Manifolds}

\indent \indent There are, at present, (at least) three inequivalent definitions of ``nearly-K\"{a}hler $6$-manifold" encountered in the literature.  We take this opportunity to compare and contrast the various notions, and also put our work in its proper context. \\

\indent  In Gray's original formulation \cite{MR0267502}, a nearly-K\"{a}hler structure on a smooth $6$-manifold $M^6$ referred to a certain kind of $\text{U}(3)$-structure on $M^6$.  A $\text{U}(3)$-structure $B \subset FM$ is equivalent to specifying on $M$ a triple $(g,J,\Omega)$ consisting of a Riemannian metric $g$, an almost-complex structure $J$, and a non-degenerate $2$-form $\Omega$ satisfying the compatibility condition $g(u,v) = \Omega(u, Jv)$.  A $6$-manifold with $\text{U}(3)$-structure is called an \textit{almost-Hermitian} $6$-manifold. \\
\indent In \cite{MR581924}, the intrinsic torsion space of a $\text{U}(3)$-structure was calculated to be of the form
$$H^{0,2}(\mathfrak{u}(3)) = \mathfrak{u}(3)^\perp \otimes \mathbb{R}^6 = W_1 \oplus W_2 \oplus W_3 \oplus W_4,$$
where $W_1, W_2, W_3, W_4$ are certain irreducible $\text{U}(3)$-modules of real dimensions $2$, $16$, $12$, $6$, respectively. \\
\indent A $\text{U}(3)$-structure was then defined to be \textit{nearly-K\"{ahler}} if its intrinsic torsion function $T \colon B \to H^{0,2}(\mathfrak{u}(3))$ takes values in $W_1$, the lowest-dimensional piece in the decomposition.  This is equivalent (see \cite{MR0267502}, \cite{MR1707644}) to requiring that $\nabla J$ satisfies $(\nabla_XJ)(X) = 0$ for all vector fields $X \in \Gamma(TM)$, or equivalently that $\nabla \Omega = \frac{1}{3}d\Omega$, where $\nabla$ is the Levi-Civita connection of the metric $g$. \\

\noindent \textit{Remark:} Note that an almost-Hermitian $6$-manifold is \textit{K\"{a}hler} if its intrinsic torsion is identically zero.  Equivalently, $\nabla J = 0$, or equivalently $\nabla \Omega = 0$. $\Box$ \\

\indent In this work, we will adopt a different definition of ``nearly-K\"{a}hler" also encountered in the literature (see, e.g., \cite{MR2253159} and \cite{MR3583352}) which entails an additional bit of structure.  For us, a ``nearly-K\"{a}hler structure" refers to a certain kind of $\text{SU}(3)$-structure. \\
\indent An $\text{SU}(3)$-structure $B \subset FM$ is equivalent to specifying on $M$ a triple $(g,J,\Omega)$ as above together with a $(3,0)$-form $\Upsilon$ such that $\Upsilon \wedge \overline{\Upsilon} = -\frac{4}{3}i\,\Omega^3$.  In fact, the data $(\Omega, \Upsilon)$, subject to appropriate algebraic conditions, is enough to reconstruct $(g,J)$.  Thus, an $\text{SU}(3)$-structure may be regarded as a pair $\Omega \in \Omega^2(M)$ and $\Upsilon \in \Omega^3(M; \mathbb{C})$ such that at each $x \in M$, there is an isomorphism $u \colon T_xM \to \mathbb{R}^6$ for which
\begin{align*}
\Omega|_x & = u^*(dx^1 \wedge dx^4 + dx^2 \wedge dx^5 + dx^3 \wedge dx^6) \\
\Upsilon|_x & = u^*(dz^1 \wedge dz^2 \wedge dz^3)
\end{align*}
where $(z^1, z^2, z^3) = (x^1 + ix^4, x^2 + ix^5, x^3 + ix^6)$ are the standard coordinates on $\mathbb{C}^3 \cong \mathbb{R}^6$. \\
\indent One can show \cite{MR1922042} that the intrinsic torsion space of an $\text{SU}(3)$-structure is of the form
\begin{align*}
H^{0,2}(\mathfrak{su}(3)) = \mathfrak{su}(3)^\perp \otimes \mathbb{R}^6 = X_0^+ \oplus X_0^- \oplus X_2^+ \oplus X_2^- \oplus X_3 \oplus X_4 \oplus X_5,
\end{align*}
where $X_0^\pm$, $X_2^\pm$, $X_3$, $X_4$, $X_5$ are certain irreducible $\text{SU}(3)$-modules of real dimensions $1$, $8$, $12$, $6$, $6$, respectively.  Following \cite{MR2287296}, we can give a more concrete description of $H^{0,2}(\mathfrak{su}(3))$ via exterior algebra.  Indeed, the $\text{SU}(3)$-modules $\Lambda^2(\mathbb{R}^6)$ and $\Lambda^3(\mathbb{R}^6)$ decompose into irreducibles as (\cite{MR2287296}, \cite{MR3656283})
\begin{align*}
\Lambda^2(\mathbb{R}^6) & = \mathbb{R}\,\Omega \oplus \Lambda^2_6 \oplus \Lambda^2_8 \\
\Lambda^3(\mathbb{R}^6) & = \mathbb{R}\,\text{Re}(\Upsilon) \oplus \mathbb{R}\,\text{Im}(\Upsilon) \oplus \Lambda^3_6 \oplus  \Lambda^3_{12},
\end{align*}
where
\begin{align*}
\Lambda^2_6 & = \{\ast(\alpha \wedge \text{Re}(\Upsilon)) \colon \alpha \in \Lambda^1\} \\
\Lambda^2_8 & = \{\beta \in \Lambda^2 \colon \beta \wedge \text{Re}(\Upsilon) = 0 \text{ and } \ast\!\beta = -\beta \wedge \Omega\} \\
\Lambda^3_6 & = \{\alpha \wedge \Omega \colon \alpha \in \Lambda^1\} = \{\gamma \in \Lambda^3 \colon \ast\! \gamma = \gamma\} \\
\Lambda^3_{12} & = \{\gamma \in \Lambda^3 \colon \gamma \wedge \Omega = 0 \text{ and } \gamma \wedge \text{Re}(\Upsilon) = 0 \text{ and } \gamma \wedge \text{Im}(\Upsilon) = 0\}.
\end{align*}
This gives the description
\begin{align*}
H^{0,2}(\mathfrak{su}(3)) \cong  \mathbb{R} \oplus \mathbb{R} \oplus \Lambda^2_8 \oplus \Lambda^2_8  \oplus \Lambda^3_{12} \oplus  \Lambda^1 \oplus  \Lambda^1.
\end{align*}
\indent It can be shown \cite{MR1922042} that the intrinsic torsion of the $\text{SU}(3)$-structure can be completely encoded in the \textit{exterior} derivatives of $\Omega$ and $\Upsilon$.  Moreover, borrowing the notation of \cite{MR3656283}, these exterior derivatives decompose as
\begin{align*}
d\Omega & = 3\tau_0\,\text{Re}(\Upsilon) + 3\hat{\tau}_0\, \text{Im}(\Upsilon) + \tau_3 + \tau_4 \wedge \Omega \\
d\,\text{Re}(\Upsilon) & = 2\hat{\tau}_0\, \Omega^2 + \tau_5 \wedge \text{Re}(\Upsilon) + \tau_2 \wedge \Omega \\
d\,\text{Im}(\Upsilon) & = -2\tau_0\, \Omega^2 - J\tau_5 \wedge \text{Re}(\Upsilon) + \hat{\tau}_2 \wedge \Omega,
\end{align*}
where $\tau_0, \hat{\tau}_0 \in \Omega^0$, $\tau_2, \hat{\tau}_2 \in \Omega^2_8$, $\tau_3 \in \Omega^3_{12}$, and $\tau_4, \tau_5 \in \Gamma(TM)$, and $\Omega^k_\ell = \Gamma(\Lambda^k_\ell(T^*M))$.  This leads to: \\

\noindent \textbf{Definition:} Let $M^6$ be a real $6$-manifold. \\
\indent A \textit{nearly-K\"{a}hler structure} on $M$ is an $\text{SU}(3)$-structure $B \subset FM$ whose intrinsic torsion function $T \colon B \to H^{0,2}(\mathfrak{su}(3))$ takes values in $X_0^-$.  \\
\indent That is: A \textit{nearly-K\"{a}hler structure} on $M$ is an $\text{SU}(3)$-structure $(\Omega, \Upsilon)$ such that $\tau_0 = \tau_2 = \hat{\tau}_2 = \tau_3 = \tau_4 = \tau_5 = 0$ and $\hat{\tau}_0 = c$ is constant.  In other words, it is an $\text{SU}(3)$-structure $(\Omega, \Upsilon)$ that satisfies
\begin{align*}
d\Omega & = 3c\,\text{Im}(\Upsilon) \\
d\,\text{Re}(\Upsilon) & = 2c\,\Omega^2 \\
d\,\text{Im}(\Upsilon) & = 0.
\end{align*}
Of course, the third equation is a consequence of the first. \\

\noindent \textit{Remark:} Note that other works (e.g. \cite{MR3583352}) instead take $\tau_0 = c$ constant and all other torsion forms equal to zero. $\Box$ \\

\indent Note that a nearly-K\"{a}hler structure has $c = 0$ if and only if it is Calabi-Yau.  Those with $c \neq 0$ are sometimes called \textit{strict} nearly-K\"{a}hler structures.  In this case, by rescaling the metric, we may take the constant $c = 1$.   For simplicity, and following \cite{MR3656283} and \cite{MR3583352}, we enact the following: \\

\noindent \textbf{Convention:} In this work, by a ``nearly-K\"{a}hler structure" \textit{we will always mean a ``strict nearly-K\"{a}hler structure, scaled so that $c = 1$."}

\subsection{The Coisotropic Orbit Condition}

\indent \indent We will need to understand how $4$-planes in $\mathbb{R}^6$ behave under the usual $\text{SU}(3)$-action.  This requires some linear algebraic preliminaries. \\

\indent Consider $(\mathbb{R}^6, g_0, \Omega_0)$ with the standard metric $g_0$, symplectic form $\Omega_0$, and orientation.  Let $\ast$ denote the corresponding Hodge star operator.  We let $(e_1, \ldots, e_6)$ be the standard basis of $\mathbb{R}^6$, and we identify $\mathbb{C}^3 \cong \mathbb{R}^6$ via $(z^1, z^2, z^3) = (x^1 + ix^4, x^2 + ix^5, x^3 + ix^6)$.  Explicitly,
\begin{align*}
g_0 & = (dx^1)^2 + \cdots + (dx^6)^2 \\
\Omega_0 & = dx^1 \wedge dx^4 + dx^2 \wedge dx^5 + dx^3 \wedge dx^6.
\end{align*}
In particular, we observe that
\begin{equation*}
\ast \Omega_0 = \textstyle \frac{1}{2}\,\Omega_0 \wedge \Omega_0. \tag{3.1}
\end{equation*}
\indent Let $V_k(\mathbb{R}^6)$ denote the Stiefel manifold of ordered orthonormal $k$-frames in $\mathbb{R}^6$, and let $\text{Gr}_k(\mathbb{R}^6)$ denote the Grassmannian of real $k$-planes in $\mathbb{R}^6$. Recall that the \textit{symplectic complement} and \textit{orthogonal complement} of a $k$-plane $E \in \text{Gr}_k(\mathbb{R}^6)$ are the respective subspaces
\begin{align*}
E^\Omega := \{v \in \mathbb{R}^6 \colon \Omega_0(v, w) = 0, \ \forall w \in E\} \\
E^\perp := \{v \in \mathbb{R}^6 \colon g_0(v, w) = 0, \ \forall w \in E\}.
\end{align*}
We say that $E$ is \textit{isotropic} if $E \subset E^\Omega$, and that $E$ is \textit{coisotropic} if $E \supset E^\Omega$.  Using (3.1), we see that for a $4$-plane $E \in \text{Gr}_4(\mathbb{R}^6)$:
$$\left.(\Omega_0 \wedge \Omega_0)\right|_E = 0 \iff E^\perp \text{ is an isotropic 2-plane} \iff E \text{ is a coisotropic 4-plane}.$$
In particular, coisotropic $4$-planes are in bijection with isotropic $2$-planes.

\pagebreak

\indent We now seek to understand the $\text{SU}(3)$-action on $2$-planes (equivalently, $4$-planes) in $\mathbb{R}^6$.  For $\theta \in [0, \pi]$, let us set
\begin{align*}
V_2(\theta) & = \text{SU}(3) \cdot \left(e_1, \cos (\theta) e_4 + \sin(\theta)e_2 \right) \subset V_2(\mathbb{R}^6) \\
\text{Gr}_2(\theta) & = \text{SU}(3) \cdot \text{span}(e_1, \cos (\theta) e_4 + \sin(\theta)e_2) \subset \text{Gr}_2(\mathbb{R}^6).
\end{align*}
Of particular interest to us is the orbit
$$\textstyle \text{Gr}_2(\frac{\pi}{2}) = \text{SU}(3) \cdot \text{span}(e_1, e_2) = \{E \in \text{Gr}_2(\mathbb{R}^6) \colon E \text{ isotropic}\}.$$ \\
\noindent \textbf{Lemma 3.1:} \\
\indent (a) Every $(v,w) \in V_2(\mathbb{R}^6)$ belongs to exactly one of the orbits $V_2(\theta)$, where $\theta \in [0,\pi]$.  The transitive $\text{SU}(3)$-actions on $V_2(0)$ and $V_2(\pi)$ have stabilizer $\text{SU}(2)$.  For $\theta \in (0,\pi)$, the transitive $\text{SU}(3)$-action on $V_2(\theta)$ is free. \\
\indent (b) Every $E \in \text{Gr}_2(\mathbb{R}^6)$ belongs to exactly one of the orbits $\text{Gr}_2(\theta)$, where $\theta \in [0,\pi)$.  The transitive $\text{SU}(3)$-action on $\text{Gr}_2(0) \cong \mathbb{CP}^2$ has stabilizer $\text{U}(2)$.  For $\theta \in (0,\pi)$, the transitive $\text{SU}(3)$-action on $\text{Gr}_2(\theta)$ has stabilizer $\text{O}(2)$. \\
\indent (c) In particular, $\text{SU}(3)$ acts transitively on
\begin{align*}
\textstyle \text{Gr}_2(\frac{\pi}{2}) \cong \{E \in \text{Gr}_4(\mathbb{R}^6) \colon E \text{ coisotropic}\} 
\end{align*}
with stabilizer
$$\text{O}(2) = \begin{pmatrix}
\cos \theta & \mp \sin \theta & 0  &  &  & \\
\sin \theta & \mp \cos \theta & 0 &  &  & \\
0 & 0 & \pm 1 & & &  \\
 &  & & \cos \theta & \pm \sin \theta & 0  \\
 &  & & \sin \theta & \mp \cos \theta & 0  \\
 &  & & 0 & 0  & \pm 1
\end{pmatrix} \leq \text{SU}(3) \leq \text{SO}(6).$$ \\
\noindent \textit{Proof:} (a) We first show that every $(v,w) \in V_2(\mathbb{R}^6)$ belongs to some $V_2(\theta)$. \\
\indent Let $(v,w) \in V_2(\mathbb{R}^6)$.  Since $\text{SU}(3)$ acts transitively on $V_1(\mathbb{R}^6) \cong \mathbb{S}^5$, there exists $A \in \text{SU}(3)$ with $Av = e_1$, so $A \cdot (v,w) = (e_1, Aw)$.  Since $Aw \perp e_1$, so $Aw \in \mathbb{R}e_4 \oplus \mathbb{C}^2$, where $\mathbb{C}^2 = \text{span}_{\mathbb{R}}(e_2, e_5, e_3, e_6)$. \\
\indent Now, the subgroup of $\text{SU}(3)$ that fixes $e_1 \in \mathbb{R}^6$ is a copy of $\text{SU}(2)$.  This $\text{SU}(2)$ acts on the orthogonal $\mathbb{R}e_4 \oplus \mathbb{C}^2$ in the usual way: it acts trivially $\mathbb{R}e_4$ and in the standard way on $\mathbb{C}^2$.  In particular, every $x \in \mathbb{R}e_4 \oplus \mathbb{C}^2$ is $\text{SU}(2)$-conjugate to an element of the form $c_4e_4 + c_2e_2$, where $c_4 \in \mathbb{R}$ and $c_2 \geq 0$. \\
\indent Thus, there exists $B \in \text{SU}(2) \leq \text{SU}(3)$ with $B \cdot Aw = c_4e_4 + c_2e_2$ for some $c_4 \in \mathbb{R}$ and $c_2 \geq 0$, so $BA \cdot (v,w) = (e_1, c_4 e_4 + c_2e_2)$.  Since $1 = \Vert w \Vert^2 = \Vert BAw \Vert^2 = c_4^2 + c_2^2$, so we may write $(c_4, c_2) = (\cos \theta, \sin \theta)$ for some $\theta \in [0, \pi]$.  Thus, $(v,w) \in V_2(\theta)$. \\
\indent To see that the orbits are disjoint, note that the composition $\Omega_0 \colon V_2(\mathbb{R}^6) \hookrightarrow \mathbb{R}^6 \times \mathbb{R}^6 \to \mathbb{R}$ is an $\text{SU}(3)$-invariant function, so is constant on the $\text{SU}(3)$-orbits $V_2(\theta)$.  Indeed,
$$\Omega_0(e_1, \cos(\theta) e_4 + \sin(\theta) e_2) = \cos(\theta).$$
In particular, if $(v,w) \in V_2(\theta_1) \cap V_2(\theta_2)$, then $\cos(\theta_1) = \cos(\theta_2)$, so $\theta_1 = \theta_2$. \\
\indent Note that $A \in \text{SU}(3)$ stabilizes $(e_1, \cos(\theta) e_4 + \sin(\theta) e_2)$ if and only if $Ae_1 = e_1$ (so $Ae_4 = e_4$) and $\sin(\theta) Ae_2 = \sin(\theta) e_2$.  For $\theta = 0$ and $\theta = \pi$, this describes $\text{SU}(2)$.  For $\theta \in (0,\pi)$, this describes the identity subgroup. \\

\indent (b) This follows from part (a) and the fibration $\text{O}(2) \to V_2(\mathbb{R}^6) \to \text{Gr}_2(\mathbb{R}^6)$. \\

\indent (c) Note that if $A \in \text{SU}(3)$ stabilizes $\text{span}(e_1, e_2)$, then $A$ also stabilizes $\text{span}(e_4, e_5)$, which forces $A$ to lie in the $\text{O}(2)$ subgroup described above. $\lozenge$ \\

\indent Thus, there are two geometrically natural first-order conditions that one could impose on the real $4$-folds in a nearly-K\"{a}hler $6$-manifold.  In one direction, we could ask that the $4$-fold be pseudo-holomorphic (normal planes lie in $\text{Gr}_2(0)$).  However, such submanifolds do not exist, even locally \cite{MR664494}.  In the other direction, we could ask that the $4$-fold be coisotropic (normal planes lie in $\text{Gr}_2(\frac{\pi}{2})$). \\
\indent There is, however, another reason to study coisotropic $4$-folds: any \textit{complete} nearly-K\"{a}hler $6$-manifold of cohomogeneity-two must have coisotropic principal orbits, as we now show. \\

\noindent \textbf{Lemma 3.2:} Let $N^n$ be a compact $G$-homogeneous Riemannian manifold.  If $\chi \in \Omega^n(N)$ is a $G$-invariant exact $n$-form on $N$,  then $\chi = 0$. \\

\noindent \textit{Proof:} Let $\chi$ be such a $G$-invariant exact $n$-form.  Write $\chi = f\,\text{vol}_N$ for some function $f \in C^\infty(N)$.  Since $\chi$ is $G$-invariant, so $f$ is $G$-invariant.  Since the $G$-action is transitive, so $f$ is constant.  Since $N$ is compact and $\chi$ is exact, Stokes' Theorem gives
$$0 = \int_N \chi = \int_N f\,\text{vol}_N = f \cdot \text{vol}(N).$$
Thus, $f = 0$, whence $\chi = 0$. $\lozenge$ \\

\noindent \textbf{Proposition 3.3:} Let $M^6$ be a nearly-K\"{a}hler $6$-manifold equipped with a $G$-action of cohomogeneity-two that preserves the $\text{SU}(3)$-structure, where $G \leq \text{Isom}(M, g)$ is closed. \\
\indent If $M$ is complete, then $M$ is compact, $G$ is compact, the quotient space $M/G$ is compact Hausdorff, and the principal $G$-orbits in $M$ are coisotropic. \\

\noindent \textit{Proof:} Suppose $M$ is complete.  Since $M$ is Einstein of positive scalar curvature, by Bonnet-Myers, $M$ is compact.  By Myers-Steenrod \cite{MR1336823}, the isometry group $\text{Isom}(M,g)$ is compact, so $G$ is compact. \\
\indent Let $N^4$ be any principal $G$-orbit in $M$.  Note that $N$ is a compact, $G$-homogeneous Riemannian manifold.  Moreover, $\Omega^2 = \frac{1}{2} d(\text{Im}(\Upsilon))$ is a $G$-invariant exact $4$-form on $N$.  Thus, by Lemma 3.2, we have $\Omega^2|_N = 0$, meaning that $N$ is coisotropic. $\lozenge$ \\

\noindent \textit{Remark:} If, moreover, $M$ is connected and simply-connected, and the Lie group $G$ is connected, then the quotient space $M/G$ is simply-connected.  See, e.g., \cite{MR0413144}. $\Box$ \\

\indent Finally, although we will not need it here, we remark that the same argument establishes: \\

\noindent \textbf{Proposition 3.4:} Let $M^6$ be a nearly-K\"{a}hler $6$-manifold equipped with a $G$-action of cohomogeneity-three that preserves the $\text{SU}(3)$-structure, where $G \leq \text{Isom}(M,g)$ is closed. \\
\indent If $M$ is complete, then $M$ is compact, $G$ is compact, the quotient space $M/G$ is compact Hausdorff, and the $3$-form $\text{Im}(\Upsilon)$ vanishes on the principal $G$-orbits.

\section{Moving Frame Setup}

\subsection{The First Structure Equations of a Nearly-K\"{a}hler $6$-Manifold}

\indent \indent Let $\pi \colon B \to M$ be an $\text{SU}(3)$-structure on a $6$-manifold $M$.  Let $\omega = (\omega^1, \ldots, \omega^6) \in \Omega^1(B; \mathbb{R}^6)$ denote the tautological $1$-form.  We will identify $\mathbb{C}^3 \cong \mathbb{R}^6$ via
$$(z^1, z^2, z^3) = (x^1 + ix^4, x^2 + ix^5, x^3 + ix^6)$$
and let $\zeta = (\zeta^1, \zeta^2, \zeta^3) \in \Omega^1(B; \mathbb{C}^3)$ denote the $\mathbb{C}$-valued tautological $1$-form:
\begin{align*}
(\zeta^1, \zeta^2, \zeta^3) & = (\omega^1 + i\omega^4, \omega^2 + i\omega^5, \omega^3 + i\omega^6).
\end{align*}
Since $B$ is an $\text{SU}(3)$-structure, the $6$-manifold $M$ is endowed with a metric $g$, a non-degenerate $2$-form $\Omega$, and a complex volume form $\Upsilon$.  Pulled up to $B$, these are exactly:
\begin{align*}
\pi^*g & = \textstyle \sum (\zeta_j \circ \overline{\zeta}_j)^2 = (\omega^1)^2 + \cdots + (\omega^6)^2 \\
\pi^*\Omega & = \textstyle \frac{i}{2} \sum \zeta_j \wedge \overline{\zeta}_j = \omega^{14} + \omega^{25} + \omega^{36} \\
\pi^*\Upsilon & = \zeta_1 \wedge \zeta_2 \wedge \zeta_3 = (\omega^1 + i\omega^4) \wedge (\omega^2 + i\omega^5) \wedge (\omega^3 + i \omega^6).
\end{align*}
\indent In the special case where the $\text{SU}(3)$-structure $B$ is nearly-K\"{a}hler, the exterior derivatives $d\zeta_i$ satisfy the \textit{first structure equations} (see \cite{MR2253159}, \cite{MR2711912}) given by
\begin{equation*}
d\zeta_i = -\kappa_{i\overline{\ell}} \wedge \zeta_\ell + \overline{\zeta_j \wedge \zeta_k}, \tag{4.1}
\end{equation*}
where $\kappa = (\kappa_{i \overline{\ell}}) \in \Omega^1(B; \mathfrak{su}(3))$ is a connection $1$-form, and where $(i,j,k)$ is an even permutation of $(1,2,3$).  In terms of the basis $(\omega^1, \ldots, \omega^6)$ for the $\pi$-semibasic $1$-forms, the structure equations (4.1) read
\begin{equation*}
d\begin{pmatrix} \omega^1 \\ \omega^2 \\ \omega^3 \\ \omega^4 \\ \omega^5 \\ \omega^6 \end{pmatrix} = -\left(\begin{array}{c c c | c c c}
0 & \alpha_3 & -\alpha_2 & -\beta_{11} & -\beta_{12} & -\beta_{13} \\
-\alpha_3 & 0 & \alpha_1 & -\beta_{21} & -\beta_{22} & -\beta_{23} \\
\alpha_2 & -\alpha_1 & 0 & -\beta_{31} & -\beta_{32} & -\beta_{33} \\ \hline
\beta_{11} & \beta_{12} & \beta_{13} & 0 & \alpha_3 & -\alpha_2 \\
\beta_{21} & \beta_{22} & \beta_{23} & -\alpha_3 & 0 & \alpha_1 \\
\beta_{31} & \beta_{32} & \beta_{33} & \alpha_2 & -\alpha_1 & 0
\end{array} \right)
\wedge \begin{pmatrix} \omega^1 \\ \omega^2 \\ \omega^3 \\ \omega^4 \\ \omega^5 \\ \omega^6 \end{pmatrix} + \begin{pmatrix} \omega^{23} - \omega^{56} \\ -\omega^{13} + \omega^{46} \\ \omega^{12} - \omega^{45} \\ -\omega^{26} + \omega^{35} \\ \omega^{16} - \omega^{34} \\ -\omega^{15} + \omega^{24} \end{pmatrix} \tag{4.2}
\end{equation*}
where $\alpha_i, \beta_{ij} \in \Omega^1(B; \mathbb{R})$ are connection $1$-forms with $\beta_{ij} = \beta_{ji}$ and $\sum \beta_{ii} = 0$. \\

\indent In this work, however, it will be convenient to express (4.1) in a different form.  Indeed, in light of the $\text{O}(2)$-representation on $\mathbb{R}^6$ described in Lemma 3.1(c), we will often prefer the basis of $\pi$-semibasic $1$-forms given by $(\eta, \overline{\eta}, \omega^3, \theta, \overline{\theta}, \omega^6)$, where
\begin{align*}
\eta & = \omega^1 + i\omega^2 \\
\theta & = \omega^4 + i\omega^5.
\end{align*}
In terms of this basis, (4.1) is equivalent to
\begin{equation*}
d\begin{pmatrix} \eta \\ \omega^3 \\ \theta \\ \omega^6 \end{pmatrix} = -\left(\begin{array}{c c c | c c c}
-i\alpha & 0 & 2i\xi_1 & -\xi_0 & -\xi_3 & -2\xi_2 \\
i\overline{\xi}_1 & -i\xi_1 & 0 & -\overline{\xi}_2 & -\xi_2 & 2\xi_0 \\ \hline
\xi_0               & \xi_3 & 2\xi_2                   &    -i\alpha & 0 & 2i\xi_1 \\
\overline{\xi}_2 & \xi_2 & -2\xi_0               & i\overline{\xi}_1 & -i\xi_1 & 0
\end{array} \right)
\wedge \begin{pmatrix}  \eta \\ \overline{\eta} \\ \omega^3 \\ \theta \\ \overline{\theta} \\ \omega^6 \end{pmatrix} +
i \begin{pmatrix}
-\eta \wedge \omega^3  + \theta \wedge \omega^6 \\
\textstyle \frac{1}{2}(\eta \wedge \overline{\eta} - \theta \wedge \overline{\theta}) \\
-\omega^3 \wedge \theta  + \eta \wedge \omega^6 \\
\textstyle \frac{1}{2}(\overline{\theta} \wedge \eta - \theta \wedge \overline{\eta}) \end{pmatrix}\!, \tag{4.3}
\end{equation*}
where $\alpha, \xi_0 \in \Omega^1(B; \mathbb{R})$ and $\xi_1, \xi_2, \xi_3 \in \Omega^1(B; \mathbb{C})$ are connection $1$-forms.  The structure equations (4.3) will be central to our calculations.

\subsection{Frame Adaptation: The $\text{O}(2)$-Bundle $P$}

\indent \indent Let $M$ be a nearly-K\"{a}hler $6$-manifold acted upon by a connected Lie group $G$ with cohomogeneity-two.  We suppose that this $G$-action preserves the $\text{SU}(3)$-structure and that the principal $G$-orbits are coisotropic.  For simplicity, the following two conventions will be in force for the rest of this work. \\

\noindent \textbf{Convention 4.1:} Without loss of generality, we suppose that $G$ acts faithfully on $M$, and that $G$ is a closed subgroup of the isometry group of $M$. $\Box$ \\

\noindent \textbf{Convention 4.2:} We restrict our attention entirely to the principal locus of $M$, by which we mean the union of principal $G$-orbits $Gx$ in $M$.  
Henceforth, when we refer to the manifold $M$, \textit{we shall always mean the principal locus of $M$}. $\Box$ \\

\indent We begin our study by adapting coframes to the foliation of $M$ by coisotropic $4$-folds.  Define the subbundle $P \subset B$ of $\text{SU}(3)$-coframes $u = (u^1, \ldots, u^6) \colon T_xM \to \mathbb{R}^6$ for which $T_xGx = \text{Ker}(u^1, u^2)$.  In other words, letting $\{e_1, \ldots, e_6\}$ denote the standard basis of $\mathbb{R}^6$, we set
$$P = \{u \in B \colon u(T_xGx) = \text{span}(e_3, e_4, e_5, e_6) \} \subset B.$$
Since $\text{SU}(3)$ acts transitively on the Grassmannian of coisotropic $4$-planes in $\mathbb{R}^6$ (Lemma 3.1), this adaptation is well-defined.  Note that $P$ is an $\text{O}(2)$-subbundle, where the inclusion $\text{O}(2) \leq \text{SU}(3) \hookrightarrow \text{GL}_6(\mathbb{R})$ is the one described in Lemma 3.1(c). \\

\noindent \textit{Remark:} The Lie group $G$ is contained in the group $\text{Aut}_{\text{O}(2)}$ of automorphisms that preserve the foliation of $M$ by coisotropic $4$-folds, which is itself contained in the full automorphism group $\text{Aut}_{\text{SU}(3)}$ of the $\text{SU}(3)$-structure:
$$G \leq \text{Aut}_{\text{O}(2)}(M) \leq \text{Aut}_{\text{SU}(3)}(M).$$
By Lemma 2.1, we see that:
\begin{align*}
4 \leq \dim(G) \leq 5, \ \ \ \ \ 4 \leq \dim(\text{Aut}_{\text{O}(2)}) \leq 7, \ \ \ \ \ 4 \leq \dim(\text{Aut}_{\text{SU}(3)}) \leq 14. \ \ \ \ \Box
\end{align*} \\
\indent Henceforth, we work on the $\text{O}(2)$-subbundle $P \subset B$ and use the same letter $\pi \colon P \to M$ to denote the restricted projection map.  Now, the connection $1$-form $\kappa \in \Omega^1(B; \mathfrak{su}(3))$ does not remain a connection form when restricted to $P$.  Indeed, for a choice of splitting $\mathfrak{su}(3) = \mathfrak{so}(2) \oplus W$, the $1$-form $\kappa|_P \in \Omega^1(P; \mathfrak{su}(3))$ decomposes as
$$\kappa|_P = \gamma_{\text{O}(2)} + \tau_{\text{O}(2)},$$
where $\gamma_{\text{O}(2)} \in \Omega^1(P; \mathfrak{so}(2))$ is a connection $1$-form and $\tau_{\text{O}(2)} \in \Omega^1(P; W)$ is $\pi$-semibasic.  In terms of the basis $(\eta, \overline{\eta}, \omega^3, \theta, \overline{\theta}, \omega^6)$, this splitting reads
\begin{equation*}
\kappa|_P = \left(\begin{array}{c c c | c c c}
-i\alpha & 0 & 0 & 0 & 0 & 0 \\
0 & 0 & 0 & 0 & 0 & 0 \\ \hline
0               & 0 & 0                   &    -i\alpha & 0 & 0 \\
0 & 0 & 0               & 0 & 0 & 0
\end{array} \right)
+
\left(\begin{array}{c c c | c c c}
0 & 0 & 2i\xi_1 & -\xi_0 & -\xi_3 & -2\xi_2 \\
i\overline{\xi}_1 & -i\xi_1 & 0 & -\overline{\xi}_2 & -\xi_2 & 2\xi_0 \\ \hline
\xi_0               & \xi_3 & 2\xi_2                   &    0 & 0 & 2i\xi_1 \\
\overline{\xi}_2 & \xi_2 & -2\xi_0               & i\overline{\xi}_1 & -i\xi_1 & 0
\end{array} \right)
\end{equation*}
In particular, $\alpha$ is a connection $1$-form for the $\text{O}(2)$-bundle $\pi \colon P \to M$, so $(\eta, \overline{\eta}, \omega^3, \theta, \overline{\theta}, \omega^6, \alpha)$ is a coframing on $P$. \\
\indent On the other hand, the $1$-forms $\xi_0 \in \Omega^1(P; \mathbb{R})$ and $\xi_1, \xi_2, \xi_3 \in \Omega^1(P; \mathbb{C})$ are $\pi$-semibasic, so that we may write
\begin{align*}
2\xi_1 & = a_{11}\eta + a_{12}\overline{\eta} + a_{13}\omega^3 + a_{14}\theta + a_{15}\overline{\theta} + a_{16}\omega^6  \\
2\xi_2 & = a_{21}\eta + a_{22}\overline{\eta} + a_{23}\omega^3  + a_{24}\theta + a_{25}\overline{\theta}+ a_{26}\omega^6 \tag{4.4} \\
\xi_3 & = a_{31}\eta + a_{32}\overline{\eta} + a_{33}\omega^3 + a_{34}\theta + a_{35}\overline{\theta} + a_{36}\omega^6 \\
\xi_0 & = a_{01}\eta + a_{02}\overline{\eta} + a_{03}\omega^3  + a_{04}\theta + a_{05}\overline{\theta} + a_{06}\omega^6
\end{align*}
for some 24 $G$-invariant functions $a_{ij} \colon P \to \mathbb{C}$.  We will refer to these $24$ functions as the \textit{torsion functions} of the $\text{O}(2)$-structure.  In the next section, we will see (Lemma 4.4) that they are not independent of one another. 

\subsection{The Torsion of the $\text{O}(2)$-Structure}

\indent \indent We continue with the setup from $\S$4.2.  The purpose of this section is to derive relations (Lemma 4.4) on the 24 functions $a_{ij} \colon P \to \mathbb{C}$ of (4.4).  In the next section, we will use this information to show (Proposition 4.5) that the acting Lie group $G$ is $4$-dimensional and non-abelian. \\

\indent We begin with the following observation: Unlike a generic $\text{O}(2)$-structure, the $\text{O}(2)$-structures in our situation enjoy a special geometric feature.  Namely, the (real) $4$-plane field $\text{Ker}(\omega^1, \omega^2) = \text{Ker}(\eta, \overline{\eta})$ on $M^6$ is integrable, its leaf space is the orbit space $\Sigma = M/G$ (recall Convention 4.2), and the quadratic form $(\omega^1)^2 + (\omega^2)^2 = \eta \circ \overline{\eta}$ descends to a Riemannian metric on $\Sigma$.  Consequently: \\

\noindent \textbf{Lemma 4.3:} There exist a $1$-form $\phi \in \Omega^1(P; \mathbb{R})$ and a function $K \in \Omega^0(P; \mathbb{R})$ such that
\begin{align*}
d\eta & = i\,\phi \wedge \eta \tag{4.5} \\
d\phi & = \textstyle \frac{i}{2}\,K\,\eta \wedge \overline{\eta}. \tag{4.6}
\end{align*} \\
\noindent \textit{Proof:} The quadratic form $\eta \circ \overline{\eta} \in \Gamma(\text{Sym}^2(T^*P))$ is both $
\text{O}(2)$-invariant and $G$-invariant, so descends to a Riemannian metric $g_\Sigma$ on the surface $\Sigma$.  Let $\varpi \colon F \to \Sigma$ denote the orthonormal frame bundle of $g_\Sigma$, and let $\widetilde{\eta} \in \Omega^1(F; \mathbb{C})$ denote the $\mathbb{C}$-valued tautological $1$-form on $F$. \\
\indent By the Fundamental Lemma of Riemannian Geometry, there exists a unique $1$-form $\widetilde{\phi} \in \Omega^1(F)$, the Levi-Civita connection of $g_\Sigma$, for which
$$d\widetilde{\eta} = i\,\widetilde{\phi} \wedge \widetilde{\eta}.$$
Now, the quotient map $\text{pr} \colon M \to \Sigma$ induces a map $\widetilde{\text{pr}} \colon P \to F$ via $\widetilde{\text{pr}}(u)(v) := u(\widetilde{v})$, where $\widetilde{v} \in TM$ is the horizontal lift of $v \in T\Sigma$.  Unwinding the definitions shows that $\widetilde{\text{pr}}^*(\widetilde{\eta}) = \eta$, whence the equation $d\eta = i\,\widetilde{\text{pr}}^*(\widetilde{\phi}) \wedge \eta$ holds on $P$.  Setting $\phi := \widetilde{\text{pr}}^*(\widetilde{\phi})$ establishes (4.5). \\
\indent Equation (4.6) now follows by differentiating (4.5).  That is, $K$ is the Gauss curvature of $g_\Sigma$. $\lozenge$ \\

\noindent \textbf{Lemma 4.4:} There exist seven $\mathbb{C}$-valued functions $p_1, p_2, p_3, p_4, q_1, q_2, r \colon P \to \mathbb{C}$ and four $\mathbb{R}$-valued functions $h_1, h_2, h_3, h_4 \colon P \to \mathbb{R}$ for which
\begin{alignat*}{6}
2\xi_1 & = & h_1\eta &{}-{} ih_2\theta & + q_1\overline{\theta} & + p_1\omega^3 &{} + p_2\omega^6 \\
2\xi_2 & = \ &{}ih_4\eta {}&+{}(h_3 - i)\,\theta &{}+ q_2\overline{\theta} &{}- ip_2\omega^3 &{} + p_3\omega^6 \tag{4.7} \\
\xi_3 & = & &\ \ \ \  p_4\theta &{} - ir\overline{\theta} & {}- iq_1\omega^3 &{} + q_2 \omega^6 \\
\xi_0 & = & &\ \ \ \ \overline{p}_4\theta &{}+ p_4\overline{\theta} &{}- h_2\omega^3 &{} + h_3\omega^6.
\end{alignat*}
Their exterior derivatives modulo $\langle \omega^1, \omega^2 \rangle = \langle \eta, \overline{\eta} \rangle$ satisfy
\begin{align*}
dp_1 & \equiv ip_1\phi & dq_1 & \equiv 2iq_1\phi & dh_1 & \equiv 0 \\
dp_2 & \equiv ip_2\phi & dq_2 & \equiv 2iq_2\phi & dh_2 & \equiv 0 \tag{4.8} \\
dp_3 & \equiv ip_3\phi &  &  & dh_3 & \equiv 0 \\
dp_4 & \equiv ip_4\phi & dr & \equiv 3ir\phi   & dh_4 & \equiv 0.
\end{align*}
Moreover, we have the formula
\begin{equation*}
\phi = \alpha + (h_1 + 1)\,\omega^3 - h_4\omega^6. \tag{4.9}
\end{equation*} \\
\indent The upshot is that we have re-expressed the torsion of the $\text{O}(2)$-structure in terms of just seven $\mathbb{C}$-valued functions and four $\mathbb{R}$-valued functions $p_i, q_i, r, h_i$ on $P$, all of which are $G$-invariant and $\text{O}(2)$-equivariant (for the $\text{O}(2)$-actions indicated by (4.8)).  Accordingly, we will refer to
$$T = (p_1, p_2, p_3, p_4, q_1, q_2, r, h_1, h_2, h_3, h_4) \colon P \to \mathbb{C}^7 \oplus \mathbb{R}^4$$
as the (\textit{intrinsic}) \textit{torsion of the $\text{O}(2)$-structure}.  Geometrically, the function $T$ describes the $1$-jet of the $\text{O}(2)$-structure (or the $2$-jet of the underlying $\text{SU}(3)$-structure) up to diffeomorphism. \\
\indent Equation (4.9) shows in particular that $(\omega^1, \ldots, \omega^6, \phi) \colon TP \to \mathbb{R}^7$ is a coframing on $P$.  Going forward, we prefer to work with the coframings $(\omega^1, \ldots, \omega^6, \phi)$ and $(\eta, \overline{\eta}, \omega^3, \theta, \overline{\theta},$ $\omega^6, \phi)$ rather than with the original $(\omega^1, \ldots, \omega^6, \alpha)$ and $(\eta, \overline{\eta}, \omega^3, \theta, \overline{\theta}, \omega^6, \alpha)$. \\

\noindent \textit{Proof of Lemma 4.4:} From (4.3) and (4.4), we have
\begin{align*}
d\eta & = i(\alpha + ia_{01} \theta + ia_{31} \overline{\theta} + (a_{11} + 1)\,\omega^3 + ia_{21}\,\omega^6) \wedge \eta \tag{4.10} \\
& - ia_{12}\,\overline{\eta} \wedge \omega^3 + a_{22}\,\overline{\eta} \wedge \omega^6 + a_{32}\,\overline{\eta} \wedge \overline{\theta} + a_{02}\,\overline{\eta} \wedge \theta    \\
& + (a_{05} - a_{34})\,\overline{\theta} \wedge \theta + (a_{03} + ia_{14})\, \omega^3 \wedge \theta  + (a_{33} + ia_{15})\,\omega^3 \wedge \overline{\theta} \\
& + (a_{24} - a_{06} + i)\,\theta \wedge \omega^6 - (a_{25} - a_{36})\,\omega^6 \wedge \overline{\theta} + (ia_{16} + a_{23})\,\omega^3 \wedge \omega^6.
\end{align*}
Equating this with $d\eta = i\phi \wedge \eta$ yields the following relations:
\begin{align*}
a_{12} & = 0 & a_{32} & = 0 & a_{15} & = ia_{33} & a_{05} & = a_{34} & a_{25} & = a_{36} \\
a_{22} & = 0 & a_{02} & = 0 & a_{06} & = a_{24} + i & a_{03} & = -ia_{14} & a_{16} & = ia_{23}
\end{align*}
From these relations, we may define
\begin{align*}
p_1 & = a_{13} & q_1 & = a_{15} = ia_{33} & h_1 & = a_{11} \\
p_2 & = a_{16} = ia_{23} & q_2 & = a_{25} = a_{36} & h_2 & = -a_{03} = ia_{14} \\
p_3 & = a_{26} &   &   & h_3 & = a_{06} = a_{24} + i \\
p_4 & = a_{05} = a_{34} & r & = ia_{35} & h_4 & = -ia_{21}.
\end{align*}
Moreover, since $\xi_0$ is real-valued, we see that $a_{03}$ and $a_{06}$ are real-valued, whence $h_2$ and $h_3$ are real-valued.  The reality of $\xi_0$ also yields $a_{01} = \overline{a}_{02} = 0$ and $a_{04} = \overline{a}_{05} = \overline{p}_4$. \\
\indent In this new notation, (4.10) and its complex conjugate now read as follows:
\begin{align*}
d\eta & = i(\alpha + ia_{31}\overline{\theta} + (h_1 + 1)\,\omega^3 - h_4\,\omega^6) \wedge \eta \\
d\overline{\eta} & = -i(\alpha - i\overline{a}_{31}\,\theta + (\overline{h}_1 + 1)\,\omega^3 - \overline{h}_4\,\omega^6) \wedge \overline{\eta}.
\end{align*}
Again equating with $d\eta = i\phi \wedge \eta$ and $d\overline{\eta} = -i\phi \wedge \overline{\eta}$, we see that $a_{31} = 0$, that $h_1$ and $h_4$ are real-valued, and that
\begin{equation*}
\phi = \alpha + (h_1 + 1)\,\omega^3 - h_4\omega^6.
\end{equation*}
This proves (4.7) and (4.9).  The proof of (4.8) is a direct calculation. $\lozenge$

\subsection{The Acting Lie Group $G$}

\noindent \textbf{Proposition 4.5:} The Lie group $G$ is $4$-dimensional and non-abelian.  In particular, if $M$ is complete, then both $G$ and the principal $G$-orbits in $M$ are finite quotients of $\text{SU}(2) \times \text{U}(1) \cong \mathbb{S}^3 \times \mathbb{S}^1$. \\

\noindent \textit{Proof:} For $X \in \mathfrak{g}$, let $X^\# \in \Gamma(TP)$ be the corresponding $G$-action vector field on $P$, by which we we mean $X^\#|_p = \left.\frac{d}{dt}\right|_{t = 0} (\exp tX) \cdot p$. \\
\indent Since $X^\#$ is tangent to the pre-images $\pi^{-1}(Gx) \subset P$, we have $\omega^1(X^\#) = \omega^2(X^\#) = 0$.  From the real part (4.5), we have $d\omega^1 = -\phi \wedge \omega^2$, whence
$$0 = \mathcal{L}_{X^\#}\omega^1  = \iota_{X^\#}(d\omega^1) + d(\iota_{X^\#}\omega^1) = \iota_{X^\#}(-\phi \wedge \omega^2) = -\phi(X^\#)\,\omega^2,$$
whence $\phi(X^\#) = 0$.   Thus, at each $p \in P$, we have
\begin{equation*}
\mathfrak{g} \cong \{X^\#|_p \in T_pP \colon X \in \mathfrak{g}\} \subset \text{Ker}\!\left.(\omega^1, \omega^2, \phi)\right|_p, \tag{4.11}
\end{equation*}
whence $\dim(G) \leq 4$.  Since $\dim(G) \geq 4$, we have equality.  In particular, the inclusion in (4.11) is an equality, so the $G$-orbits in $P$ are the integral $4$-folds of $\mathcal{I}_G := \langle \omega^1, \omega^2, \phi \rangle$. \\
\indent Let us now identify $G$ with an integral $4$-fold of $\mathcal{I}_G$.  Via this identification, $\{\omega^3, \omega^4, \omega^5$, $\omega^6\}$ is a basis of left-invariant $1$-forms on $G$.  Let $\{X_3, X_4, X_5, X_6\}$ be a basis of $\mathfrak{g} = \{\text{left-invariant vector}$ $\text{fields on }G\}$ whose dual basis is $\{\omega^3, \omega^4, \omega^5, \omega^6\}$.  From (4.3) and (4.7), one may calculate that
\begin{align*}
d\omega^3 & \equiv \textstyle -\frac{3}{2}i\,\theta \wedge \overline{\theta} = -3\,\omega^4 \wedge \omega^5 \ \ \ (\text{mod } \eta, \overline{\eta}, \omega^3, \omega^6, \phi).
\end{align*}
Thus, $[X_4, X_5] = 3X_3$, so $G$ is non-abelian. \\
\indent If $M$ is complete, then (Proposition 3.3) $G$ is a compact $4$-dimensional non-abelian Lie group.  Hence, by the classification of compact Lie groups \cite{MR0413144}, $G$ must be a finite quotient of $\text{SU}(2) \times \text{U}(1)$.  In this case, since the principal $G$-orbits are $4$-dimensional $G$-homogeneous spaces, they must also be finite quotients of $\text{SU}(2) \times \text{U}(1)$. $\lozenge$

\subsection{Geometric Interpretation of the Torsion}

\indent \indent We pause to interpret the torsion functions $p_1, p_2, p_3, p_4, q_1, q_2, r, h_1, h_2, h_3, h_4$ geometrically.  This section is parenthetical to the rest of this work: the main results in $\S$5 will not draw on these remarks.

\subsubsection{Background: Riemannian Submersions}

\indent \indent Let $\text{pr} \colon (M^n, g) \to (\Sigma^k, g_\Sigma)$ denote an arbitrary Riemannian submersion between Riemannian manifolds $M$ and $\Sigma$.  Recall that the vertical distribution $\mathcal{V} = \text{Ker}(\text{pr}_*)$ is given by the tangent spaces to the $\text{pr}$-fibers, and the horizontal distribution $\mathcal{H} = \mathcal{V}^\perp$ is the orthogonal complement. \\
\indent The geometry of the submersion $\text{pr}$ is governed by the two $(1,2)$-tensor fields on $M$, called the $\mathsf{A}$-tensor and $\mathsf{T}$-tensor, given by
\begin{align*}
\mathsf{A}(X,Y) & = (\nabla_{X^{\text{Hor}} } Y^{\text{Ver}})^{\text{Hor}} + (\nabla_{X^{\text{Hor}} } Y^{\text{Hor}})^{\text{Ver}} \\
\mathsf{T}(X,Y) & = (\nabla_{X^{\text{Ver}} } Y^{\text{Ver}})^{\text{Hor}} + (\nabla_{X^{\text{Ver}} } Y^{\text{Hor}})^{\text{Ver}}
\end{align*}
for $X,Y \in \Gamma(TM)$, where $\nabla$ is the Levi-Civita connection on $TM$, and where $\text{Ver} \colon TM \to \mathcal{V}$ and $\text{Hor} \colon TM \to \mathcal{H}$ are the projections onto the vertical and horizontal distributions, respectively. \\
\indent Note that $\mathsf{A} \equiv 0$ if and only if the horizontal distribution $\mathcal{H}$ is integrable.  Indeed, $\mathsf{A}(X,Y) = \frac{1}{2}[X,Y]^{\text{Ver}}$ for $X,Y \in \Gamma(\mathcal{H})$.  Meanwhile, the $\mathsf{T}$-tensor is essentially the second fundamental form $\text{I\!I}$ of the $\text{pr}$-fibers.  Indeed, $\mathsf{T}(X,Y) = \text{I\!I}(X,Y)$ for $X,Y \in \Gamma(\mathcal{V})$. \\
\indent Finally, we point out that $\mathsf{A}(X, \cdot) = 0$ for all $X \in \Gamma(\mathcal{V})$, and similarly $\mathsf{T}(X, \cdot) = 0$ for all $X \in \Gamma(\mathcal{H})$.  Thus, the $\mathsf{A}$- and $\mathsf{T}$-tensors are recovered, respectively, from the knowledge of $\mathsf{A}(X, \cdot)$ for $X \in \Gamma(\mathcal{H})$ and $\mathsf{T}(X, \cdot)$ for $X \in \Gamma(\mathcal{V})$.  For more information, see \cite{MR2371700}.

\subsubsection{Geometric Interpretation of the Torsion Functions}

\indent \indent We now return to our usual setting, which is that of a cohomogeneity-two nearly-K\"{a}hler $6$-manifold $(M^6,g)$ with coisotropic principal orbits. \\
\indent Let $\text{pr} \colon (M^6, g) \to \Sigma$ denote the projection to the orbit space.  As in Lemma 4.3, we equip $\Sigma$ with the Riemannian metric $g_\Sigma = \eta \circ \overline{\eta}$, so that $\text{pr}$ is a Riemannian submersion.  We claim that: \\

\noindent \textbf{Proposition 4.6:} \\
\indent (a) The torsion functions $h_1, h_4$ determine the $\mathsf{A}$-tensor, and conversely. \\
\indent (b) The torsion functions $p_1, p_2, p_3, p_4, q_1, q_2, r, h_2, h_3$ determine the $\mathsf{T}$-tensor, and conversely. \\

\indent To see this, let $F_{\text{SO}(6)}M$ denote the oriented orthonormal frame bundle of the metric $g$.  By the Fundamental Lemma of Riemannian Geometry, there is a unique $1$-form $\psi \in \Omega^1(F_{\text{SO}(6)}M; \mathfrak{so}(6))$, the Levi-Civita connection, such that
$$d\omega = -\psi \wedge \omega.$$
One can check that, using the notation of (4.2), the Levi-Civita connection (restricted to $P$) is
$$\psi = \left(\begin{array}{c c | c c c c}
0 & \alpha_3 + \frac{1}{2}\omega^3 & -\alpha_2 - \frac{1}{2}\omega^2 & -\beta_{11} & -\beta_{12} - \frac{1}{2}\omega^6 & -\beta_{13} + \frac{1}{2}\omega^5  \\
-\alpha_3 - \frac{1}{2}\omega^3  & 0 & \alpha_1 + \frac{1}{2}\omega^1 & -\beta_{12} + \frac{1}{2}\omega^6 & -\beta_{22} & -\beta_{23} - \frac{1}{2}\omega^4  \\ \hline
\alpha_2 + \frac{1}{2}\omega^2 & -\alpha_1 - \frac{1}{2}\omega^1 & 0 & -\beta_{13} - \frac{1}{2}\omega^5 & -\beta_{23} + \frac{1}{2}\omega^4 & -\beta_{33} \\
\beta_{11} & \beta_{12} - \frac{1}{2}\omega^6 & \beta_{13} + \frac{1}{2}\omega^5 & 0 & \alpha_3 - \frac{1}{2}\omega^3 & -\alpha_2 + \frac{1}{2}\omega^2 \\
\beta_{12} + \frac{1}{2}\omega^6 & \beta_{22} & \beta_{23} - \frac{1}{2}\omega^4 & -\alpha_3 + \frac{1}{2}\omega^3 & 0 & \alpha_1 - \frac{1}{2}\omega^1 \\
\beta_{13} - \frac{1}{2}\omega^5 & \beta_{23} - \frac{1}{2}\omega^4 & \beta_{33} & \alpha_2 - \frac{1}{2}\omega^2 & -\alpha_1 + \frac{1}{2}\omega^1 & 0
\end{array} \right)\!.$$
Thus, letting $\nabla$ denote the corresponding covariant derivative operator on $TM$, we have
\begin{align*}
(\nabla_X e_i)^{ \text{Hor} } & = \psi^1_i(X)e_1 + \psi^2_i(X)e_2 \\
(\nabla_X e_i)^{ \text{Ver} } & = \psi^3_i(X)e_3 + \psi^4_i(X)e_4 + \psi^5_i(X)e_5 + \psi^6_i(X)e_6
\end{align*}
where $(e_1, \ldots, e_6)$ is any local $\text{O}(2)$-frame field.  Using these formulas, together with (4.7), one can compute the $\mathsf{A}$- and $\mathsf{T}$-tensors. \\
\indent For example, a calculation shows that
\begin{align*}
\mathsf{A}(e_1 + ie_2, e_1) & = \textstyle i\left[ (h_1 + \frac{1}{2})e_3 - h_4 e_6 \right] & \mathsf{A}(e_1 + ie_2, e_4 ) & = 0 \\
\mathsf{A}(e_1 + ie_2, e_2) & = \textstyle -i\left[ (h_1 + \frac{1}{2})e_3 - h_4 e_6 \right] & \mathsf{A}(e_1 + ie_2, e_5 ) & = 0 \\
\mathsf{A}(e_1 + ie_2, e_3) & = \textstyle -i(h_1 + \frac{1}{2})\,(e_1 + ie_2) & \mathsf{A}(e_1 + ie_2, e_6) & = ih_4\,(e_1 + ie_2)
\end{align*}
where we have extended $\mathsf{A}$ to be $\mathbb{C}$-bilinear.  In particular, we observe that the $2$-plane field $\text{Ker}(\omega^3, \omega^4, \omega^5, \omega^6)$ in $M$ normal to the principal $G$-orbits is integrable if and only if $\mathsf{A} = 0$, or equivalently
$$\textstyle h_1 = -\frac{1}{2} \ \ \ \text{ and } \ \ \ h_4 = 0.$$
Computations by the author suggest that this integrability cannot happen (locally), but the details require closer examination. \\

\indent We now exhibit the second fundamental form $\text{I\!I}$ of the principal orbits.  This is the normal bundle-valued quadratic form $\text{I\!I} \in \Gamma(\text{Sym}^2(T^*M) \otimes NM)$ given by
$$\text{I\!I} = h^1_{ij}\,\omega^i \omega^j \otimes e_1 + h^2_{ij}\, \omega^i \omega^j \otimes e_2$$
where
\begin{align*}
\left(h^1_{ij}\right) & = \begin{bmatrix}
-\text{Im}(p_1) & -\text{Im}(q_1) + h_2 & \text{Re}(q_1) & -\text{Im}(p_2) \\
-\text{Im}(q_1) + h_2 & -3\,\text{Re}(p_4) - \text{Im}(r) & -\text{Im}(p_4) + \text{Re}(r) & -\text{Re}(q_2) - h_3 \\
\text{Re}(q_1) & -\text{Im}(p_4) + \text{Re}(r) & -\text{Re}(p_4) + \text{Im}(r) & -\text{Im}(q_2) - \frac{1}{2} \\
-\text{Im}(p_2) & -\text{Re}(q_2) - h_3 & -\text{Im}(q_2) - \frac{1}{2} & -\text{Re}(p_3)
\end{bmatrix} \\
\\
\left(h^2_{ij}\right) & = \begin{bmatrix}
\text{Re}(p_1) & \text{Re}(q_1) & \text{Im}(q_1) + h_2 & \text{Re}(p_2) \\
\text{Re}(q_1) & -\text{Im}(p_4) + \text{Re}(r) & -\text{Re}(p_4) + \text{Im}(r) & -\text{Im}(q_2) + \frac{1}{2} \\
\text{Im}(q_1) + h_2 & -\text{Re}(p_4) + \text{Im}(r) & -3\,\text{Im}(p_4) - \text{Re}(r) & \text{Re}(q_2) - h_3 \\
\text{Re}(p_2) & -\text{Im}(q_2) + \frac{1}{2} & \text{Re}(q_2) - h_3 & -\text{Im}(p_3)
\end{bmatrix}\!.
\end{align*}
Conversely, one can invert these formulas to recover $p_1, p_2, p_3, p_4, q_1, q_2, r, h_2, h_3$ in terms of $h^1_{ij}, h^2_{ij}$.  Indeed,
\begin{align*}
p_1 & = h^2_{11} - ih^1_{11} & 2q_1 & = 2h^1_{13} + i(h^2_{13} - h^1_{12}) & 2h_2 & = h^1_{12} + h^2_{13} \\
p_2 & = h^2_{14} - ih^1_{14} & 2q_2 & = (h^2_{34} - h^1_{24}) - i(1 + 2h^1_{34}) & 2h_3 & = -h^1_{24} - h^2_{34} \\ 
p_3 & = -h^1_{44} - ih^2_{44} &  \\
4p_4 & = -(h^1_{22} + h^1_{33}) - i(h^1_{23} + h^2_{33})  & 4r & = (3h^1_{23} - h^2_{33}) + i(3h^1_{33} - h^1_{22})
\end{align*}
illustrating that these torsion functions are simply affine-linear combinations of second fundamental form coefficients.  We also note that the mean curvature of the principal orbits is
\begin{equation*}
H = -( \text{Im}(p_1) + \text{Re}(p_3) + 4\,\text{Re}(p_4)  )\,e_1 + ( \text{Re}(p_1) - \text{Im}(p_3) - 4\,\text{Im}(p_4) )\,e_2,
\end{equation*}
so that a principal orbit is minimal if and only if $p_1 + ip_3 + 4ip_4 = 0$.  Finally, we point out that the presence of the $\frac{1}{2}$ terms in $(h^1_{ij})$ and $(h^2_{ij})$ implies that they cannot vanish simultaneously, so that none of the principal orbits can be totally-geodesic.

\subsubsection{Descent to $M$}

\indent \indent We caution the reader that the $1$-forms $\eta, \overline{\eta}, \omega^3, \theta, \overline{\theta}, \omega^6$ and torsion functions $p_j$, $q_j$, $r$, $h_j$ are defined on the bundle $P$, not on the base manifold $M$.  However, $\text{O}(2)$-invariant combinations of these will descend to be well-defined (possibly up to sign) on $M$. \\
\indent For example, the quadratic forms $\eta \circ \overline{\eta}$, $(\omega^3)^2$, $\theta \circ \overline{\theta}$, and $(\omega^6)^2$ descend to $M$, and the differential forms $\eta \wedge \overline{\eta}$, $\omega^3$, $\theta \wedge \overline{\theta}$, and $\omega^6$ descend to be well-defined \textit{up to sign}.
Similarly, the norms of the torsion functions $|p_1|, |p_2|, |p_3|, |p_4|$, $|q_1|, |q_2|$, $|r|$, $|h_1|, |h_2|, |h_3|, |h_4|$ are well-defined on $M$, while the $1$-forms, quadratic forms, and cubic forms
$$\overline{p}_1\eta, \ \ \ \ \ \overline{p}_2\eta,  \ \ \ \ \ \overline{p}_3 \eta,  \ \ \ \ \ \overline{p}_4 \eta, \ \ \ \ \ \overline{q}_1\,\eta \circ \eta, \ \ \ \ \ \overline{q}_2\,\eta \circ \eta, \ \ \ \ \ \overline{r} \,\eta \circ \eta \circ \eta$$
descend to be well-defined up to sign.

\subsection{Calculus on the Orbit Space $\Sigma$}

\indent \indent In this section, we describe the geometric structure on the orbit space $\Sigma$ of the principal locus in terms of our moving frame apparatus.  This will be used in $\S$5 to derive the elliptic PDE systems on $\Sigma$ satisfied by cohomogeneity-two nearly-K\"{a}hler structures.

\subsubsection{Operators on Holomorphic Line Bundles over $\Sigma$}

\indent \indent Recall that the quadratic form $\eta \circ \overline{\eta} = (\omega^1)^2 + (\omega^2)^2$ descends to a well-defined Riemannian metric $g_\Sigma$ on the quotient surface $\Sigma$.  Let $\varpi \colon F \to \Sigma$ denote the corresponding orthonormal coframe bundle over $\Sigma$.  Recall also from the proof of Lemma 4.3 that the Levi-Civita connection $\phi \in \Omega^1(F)$ and Gauss curvature $K \colon F \to \mathbb{R}$ of the metric $g_\Sigma$ satisfy
\begin{align*}
d\eta & = i\,\phi \wedge \eta \\
d\phi & = \textstyle \frac{i}{2}\,K\,\eta \wedge \overline{\eta},
\end{align*}
\noindent \textit{Remark:} Here, it would perhaps be more proper to write $\widetilde{\phi}$ and $\widetilde{\eta}$ in place of $\phi$ and $\eta$, as in the proof of Lemma 4.3.  However, we will often follow a convention common in the method of moving frames, denoting forms on $F$ and their pullbacks by $\widetilde{\text{pr}} \colon P \to F$ by the same notation. $\Box$ \\

\indent Define an almost-complex structure on $\Sigma$ as follows: For $\sigma \in \Lambda^1(\Sigma; \mathbb{C}) = T^*\Sigma^\mathbb{C}$, we declare that $\sigma \in \Lambda^{1,0}(\Sigma)$ iff $\varpi^*(\sigma) \in \text{span}_{\mathbb{C}}(\eta)$.  For dimension reasons, this almost-complex structure is integrable, and by construction, it is compatible with the metric.  Since the associated $2$-form $\frac{i}{2}\eta \wedge \overline{\eta}$ is closed (again by dimension reasons), so $\Sigma$ is K\"{a}hler.
\begin{diagram}
T^*F^{\mathbb{C}} & \rTo & T^*\Sigma^{\mathbb{C}} & \ \ \ \ \ \ \ \ \ \ \ \ \ \ & \varpi^*(K_\Sigma^n \otimes T^*\Sigma^\mathbb{C}) & \rTo & K_\Sigma^n \otimes T^*\Sigma^\mathbb{C} \\
\dTo & & \dTo&  & \dTo & & \dTo \\
F & \rTo^\varpi & \Sigma & & F & \rTo^\varpi & \Sigma
\end{diagram}
\indent By construction, $K_\Sigma = \Lambda^{1,0}(\Sigma) \to \Sigma$ is a holomorphic line bundle, as are its tensor powers $K_\Sigma^n \cong \text{Sym}^n(\Lambda^{1,0}(\Sigma))$.  In particular, each $K_\Sigma^n$ admits a $\overline{\partial}$-operator:
$$\overline{\partial} \colon \Gamma(K_\Sigma^n) \to \Gamma(K_\Sigma^n \otimes \Lambda^{0,1}(\Sigma)).$$
\indent The Levi-Civita connection $\phi \in \Gamma(T^*F^\mathbb{C})$ of the metric on $\Sigma$ induces a covariant derivative operator on $T^*\Sigma^\mathbb{C}$.  Since $\Sigma$ is K\"{a}hler, there is an induced covariant derivative operator on $K_\Sigma = \Lambda^{1,0}(\Sigma)$, and hence also on all of its tensor powers:
$$\nabla \colon \Gamma(K_\Sigma^n) \to \Gamma(K_\Sigma^n \otimes T^*\Sigma^\mathbb{C})$$
Let us give a more explicit description of $\nabla$. Let $\sigma \in \Gamma(K_\Sigma^n)$ be a smooth section, say $\varpi^*(\sigma) = f\eta^n$ for some function $f \in \Omega^0(F; \mathbb{C})$.  Write
$$df = f'\eta + f''\overline{\eta} - f_0i\phi.$$
Then $\nabla \sigma \in \Gamma(\Sigma; K^n_\Sigma \otimes T^*\Sigma^{\mathbb{C}} )$ is the section such that $\varpi^*(\nabla \sigma) \in \Gamma(F; \varpi^*(K^n_\Sigma \otimes T^*\Sigma^{\mathbb{C}} ) )$ is given by
$$\varpi^*(\nabla\sigma) = \eta^n \otimes (df + f_0i\phi) = \eta^n \otimes (f'\eta + f''\overline{\eta}).$$
\indent Since $\Sigma$ is K\"{a}hler, the Levi-Civita connection on $T^*\Sigma$ coincides with the Chern connection on $\Lambda^{1,0}\Sigma$ (under the isomorphism $T^*\Sigma \cong \Lambda^{1,0}(\Sigma)$), so that $\nabla$ is compatible with both the holomorphic structure and Hermitian structure on $K_\Sigma = \Lambda^{1,0}(\Sigma)$.  In particular,
$$\varpi^*(\overline{\partial}\sigma) = \varpi^*(\nabla^{0,1}\sigma) = f'' \eta^n \otimes \overline{\eta}.$$
\subsubsection{The Laplacian of the Conformal Factor}

\indent \indent Let $(z)$ be a local holomorphic coordinate on $\Sigma$.  We can write
$$\lambda \eta = \varpi^*(dz)$$
for some non-zero function $\lambda = e^{u+iv} \colon F \to \mathbb{C}$, where here $u,v \colon F \to \mathbb{R}$.  Note that $|\lambda|^2 = e^{2u}$, and hence $u = \frac{1}{2} \log |\lambda|^2$, both descend to well-defined functions on $\Sigma$.  We also have
\begin{align*}
e^{2u}\,\eta \circ \overline{\eta} & = \varpi^*\!\left(dz \circ d\overline{z}\right) \\
e^{2u}\,\eta \wedge \overline{\eta} & = \varpi^*\!\left(dz \wedge d\overline{z}\right)\!.
\end{align*}
\indent Let $\ast$ be the Hodge star of the metric and orientation on $\Sigma$.  The Hodge Laplacian $\Delta$ on functions is then $\Delta = -\ast d \ast d$.  The following result is classical, but we prove it for completeness. \\

\noindent \textbf{Lemma 4.7:} We have $\Delta u = -K$. \\

\noindent \textit{Proof:} We begin by observing that
\begin{align*}
i\phi \wedge \eta = d\eta = d\!\left(\frac{1}{\lambda}\,\varpi^*(dz)\right) = -\frac{d\lambda}{\lambda^2} \wedge \varpi^*(dz) & = -\frac{d\lambda}{\lambda} \wedge \eta,
\end{align*}
so
$$\left(i\phi + \frac{d\lambda}{\lambda} \right) \wedge \eta = 0.$$
Thus, by Cartan's Lemma, there exists a function $h \colon F \to \mathbb{C}$ with
$$\frac{d\lambda}{\lambda} = -i\phi + h\eta.$$
\indent Write the exterior derivative of $h$ in the form $dh = h'\eta + h'' \overline{\eta} + h_0\phi$ for some functions $h', h'', h_0 \colon F \to \mathbb{C}$.  Then
\begin{align*}
-\frac{1}{2}K\, \eta \wedge \overline{\eta} = d(i \phi) = d\!\left( - \frac{d\lambda}{\lambda} \right) + d(h\eta)  & = dh \wedge \eta + h\,d\eta = -h'' \eta \wedge \overline{\eta} + (h_0 + ih) \,\phi \wedge \eta
\end{align*}
shows that
$$dh = h'\eta + \frac{1}{2}K\overline{\eta} - ih\phi.$$
Finally, noting that
$$d(\log |\lambda|^2) = \frac{d\lambda}{\lambda} + \frac{d\overline{\lambda}}{\overline{\lambda}} = h \eta + \overline{h}\overline{\eta}$$
we may calculate
\begin{align*}
\Delta u = \frac{1}{2} \Delta\!\left(\log |\lambda|^2\right) = -\frac{1}{2}\ast d \ast d(\log |\lambda|^2) = -\frac{1}{2} \ast d \ast ( h \eta + \overline{h}\overline{\eta} ) & = -\frac{1}{2} \ast\!(iK\, \eta \wedge \overline{\eta}) = -K. 
\end{align*}
\noindent $\lozenge$

\section{Local Existence and Generality}

\indent \indent We continue with the setup of $\S$4, which we reiterate for clarity.  We let $M$ be a nearly-K\"{a}hler $6$-manifold acted upon by a connected Lie group $G$ with cohomogeneity-two.  We suppose that this $G$-action preserves the $\text{SU}(3)$-structure $(J, \Omega, \Upsilon)$ and that the principal $G$-orbits are coisotropic.  Conventions 4.1 and 4.2 (stated in $\S$4.2) remain in force. \\
\indent We continue to work on the principal $\text{O}(2)$-bundle $\pi \colon P \to M$, defined in $\S$4.2 as a frame adaptation.  On $P$, we work with either of the global coframings $(\omega^1, \ldots, \omega^6, \phi)$ or $(\eta, \overline{\eta}, \omega^3, \theta, \overline{\theta}, \omega^6, \phi)$, and their exterior derivatives are given by (4.3) and (4.7).  Finally, the intrinsic torsion of the $\text{O}(2)$-structure has been encoded as a function
$$T = (p_1, p_2, p_3, p_4, q_1, q_2, r, h_1, h_2, h_3, h_4) \colon P \to \mathbb{C}^7 \oplus \mathbb{R}^4$$
which satisfies the $\text{O}(2)$-equivariance described in (4.8). \\

\indent Our primary objective is to prove a local existence/generality theorem for nearly-K\"{a}hler $6$-manifolds of cohomogeneity-two (always assuming the principal orbits are coisotropic) by appealing to Cartan's Third Theorem (Theorem 2.2).  Concretely, this means satisfying the integrability conditions
\begin{align*}
& d(d\eta) = d(d\overline{\eta}) = 0 \\
& d(d\omega^3) = d(d\theta) = d(d\overline{\theta}) = d(d\omega^6) = 0 \tag{5.1} \\
& d(d\phi) = 0 \tag{5.2} \\
& d(dp_i) = d(dq_i) = d(dh_i) = 0 \tag{5.3} \\
& d(dr) = 0, \tag{5.4}
\end{align*}
as well as ensuring the involutivity and correct dimension of the tableau of free derivatives.  Fortunately, the equations $d(d\eta) = d(d\overline{\eta}) = 0$ are already satisfied (by Lemma 4.3). \\
\indent By contrast, the integrability conditions (5.1) are quite complicated, consisting of $4 \binom{6}{3} = 80$ quadratic equations on $55$ real-valued functions: the $18$ real and imaginary parts of the torsion functions, their $36$ ``directional derivatives" in the two directions normal to the $G$-orbits, and the Gauss curvature $K$ of the orbit space $\Sigma$.  Thus, arranging for (5.1) will occupy us for some time.

\subsubsection{The Three Types}

\indent \indent We begin by solving two of the simpler quadratic equations arising in (5.1).  Namely, we calculate
\begin{align*}
0 = d(d\theta) \wedge \eta \wedge \overline{\eta} \wedge \theta & = -8\left[ q_1(h_1 + 3 + ih_3) - q_2(h_2 + ih_4) \right] \omega^{123456} \\
0 = d(d\theta) \wedge \eta \wedge \overline{\eta} \wedge \theta & = -4i\left[ q_1\overline{q}_2 + \overline{q}_1q_2 - 2\,(h_1h_2 + 3h_2 + h_3h_4) \right] \omega^{123456}
\end{align*}
yielding the quadratics
\begin{align*}
q_1(h_1 + 3 + ih_3) - q_2(h_2 + ih_4) & = 0 \tag{5.5a} \\
q_1\overline{q}_2 + \overline{q}_1q_2 & = 2\,(h_1h_2 + 3h_2 + h_3h_4). \tag{5.5b}
\end{align*}
To solve this system, we introduce the $\mathbb{C}$-valued functions
\begin{align*}
s_1 & = (h_1 + 3) + ih_3 \\
s_2 & = h_2 + ih_4.
\end{align*}
For $z,w \in \mathbb{C}$, we let $\langle z,w \rangle = \text{Re}(z\overline{w})$ denote the euclidean inner product on $\mathbb{C} \simeq \mathbb{R}^2$, and let $\Vert z \Vert = \sqrt{z\overline{z}}$ denote the euclidean norm.   Then equations (5.5a)-(5.5b) are simply
\begin{align*}
q_1s_1 - q_2s_2 & = 0 \tag{5.6a} \\
\langle q_1, q_2 \rangle & = \langle s_2, s_1 \rangle. \tag{5.6b}
\end{align*}
The solution to (5.6a)-(5.6b) is provided by the following geometric fact. \\

\noindent \textbf{Lemma 5.1:} Let $a,b,c,d \in \mathbb{C}$ be complex numbers satisfying both 
\begin{align*}
ad - bc & = 0 \\
\langle a,b \rangle & = \langle c,d \rangle.
\end{align*}
Then exactly one of the following holds: \\
\indent (i) $a = b = c = d = 0$. \\
\indent (ii) $\langle a,b \rangle = \langle c,d \rangle \neq 0$ and $\Vert a \Vert = \Vert c \Vert$ and $\Vert b \Vert = \Vert d \Vert$.  \\
\indent (iii) $\langle a,b \rangle = \langle c,d \rangle = 0$ and $(a,b,c,d) \neq (0,0,0,0)$. \\

\indent Accordingly, we partition the class of manifolds under consideration into three types: \\

\noindent \textbf{Definition:} Let $M$ be a nearly-K\"{a}hler $6$-manifold of cohomogeneity-two with coisotropic principal orbits.  We say that $M$ is of: \\
\indent $\bullet$ \textit{Type I} if $q_1 = q_2 = s_1 = s_2 = 0$. \\
\indent $\bullet$  \textit{Type II} if $\langle q_1, q_2 \rangle = \langle s_2, s_1 \rangle \neq 0$ and $\Vert q_1 \Vert = \Vert s_2 \Vert$ and $\Vert q_2 \Vert = \Vert s_1 \Vert$. \\
\indent $\bullet$ \textit{Type III} if $\langle q_1, q_2 \rangle = \langle s_2, s_1 \rangle = 0$ and $(q_1, q_2, s_1, s_2) \neq (0,0,0,0)$. \\

\noindent \textit{Remark:} The Type conditions are pointwise: to be precise, we should speak of being ``Type I at $p \in P$," and so on.  Since the Type conditions are $\text{O}(2)$-invariant, it makes sense to speak of $M$ as being ``Type I at $m \in M$," etc.  It is conceivable for a nearly-K\"{a}hler structure on $M$ to be of (say) Type I at some points of $M$ and be of Type II at others. $\Box$ \\

\indent In the sequel, we study each Type of cohomogeneity-two nearly-K\"{a}hler structure separately.  In each case, the primary challenge will be solving the quadratic equations (5.1).  Once this is done, we will solve (5.2), (5.3), (5.4) and draw conclusions. \\
\indent We will see that the algebra involved in solving (5.1) is fairly simple for Type I structures, but is significantly more labor intensive in the Type II case, and even more so in the Type III case.

\subsection{Type I}

\indent \indent In this section, we study nearly-K\"{a}hler structures of Type I.  In particular, we prove a local existence/generality result (Theorem 5.4) for these structures.  We then show that for this Type, the acting Lie group $G$ is nilpotent (Proposition 5.5), and hence the underlying metrics are incomplete.  We conclude the section by deriving an elliptic PDE system (Proposition 5.6) on the quotient surface $\Sigma$ which Type I structures must satisfy.

\subsubsection{The Integrability Conditions}

\indent \indent Our first task is to make explicit the integrability conditions (5.1), which amount to quadratic equations on both the torsion functions $p_i, q_i, r, h_i$ and on their first derivatives.  In preparation, we express the exterior derivatives of $p_4$ and $r$ as
\begin{align*}
dp_4 & = p_4'\,\eta + p_4'' \,\overline{\eta} + ip_4\phi \\
dr & = r'\,\eta + r'' \,\overline{\eta} + 3ir\phi
\end{align*}
for some functions $p_4', p_4'', r', r'' \colon P \to \mathbb{C}$.  In doing so, we have used the $G$-invariance and $\text{O}(2)$-equivariance (4.8) of $p_4$ and $r$. \\

\noindent \textbf{Lemma 5.2:} Let $M$ be a nearly-K\"{a}hler manifold of Type I. \\
\indent (a) On the $\text{O}(2)$-coframe bundle $P$, the following algebraic relations hold:
\begin{align*}
(p_1, p_2, p_3, p_4) & = (-4ip_4, 0, 4p_4, p_4) \tag{5.7} \\
(q_1, q_2) & = (0,0) \\
(h_1, h_2, h_3, h_4) & = (-3, 0, 0, 0).
\end{align*}
Thus, the torsion can be expressed in terms of the functions $p_4$ and $r$.  \\
\indent (b) On the $\text{O}(2)$-coframe bundle $P$, the following differential relations hold:
\begin{align*}
p_4' & = \textstyle -6|p_4|^2 - \frac{3}{2} & r' & = -5ip_4^2 - \overline{p}_4r \tag{5.8} \\
p_4'' & = - 5p_4^2 + i\overline{p}_4r  & K & = 2(6 + |r|^2 - |p_4|^2).
\end{align*} \\
\noindent \textit{Proof:} The equations $(q_1, q_2) = (0,0)$ and $(h_1, h_2, h_3, h_4) = (-3, 0, 0, 0)$ are immediate from the definition of ``Type I."  For the others, we calculate (using (4.3) and (4.7))
\begin{align*}
0 = d(d\theta) \wedge \theta \wedge \overline{\theta} \wedge \omega^6 & = 8p_2\,\omega^{123456} \\
0 = d(d\theta) \wedge \eta \wedge \overline{\theta} \wedge \omega^6 & = 6(p_1 + ip_3)\,\omega^{123456} \\
0 = d(d\overline{\theta}) \wedge \eta \wedge \theta \wedge \omega^6 & = 8(p_1 + 4ip_4)\,\omega^{123456}
\end{align*}
and
\begin{align*}
0 = d(d\omega^3) \wedge \overline{\eta} \wedge \omega^3 \wedge \theta & = \textstyle 16\!\left(p_4' + 6 |p_4|^2 + \frac{3}{2}\right) \omega^{123456} \\
0 = d(d\omega^3) \wedge \eta \wedge \theta \wedge \omega^3 & = 16\!\left(p_4'' + 5p_4^2 - i\overline{p}_4r\right) \omega^{123456} \\
0 = d(d\theta) \wedge \omega^3 \wedge \overline{\theta} \wedge \omega^6 & = 2\!\left(12 + 2 |r|^2 - 2|p_4|^2 - K \right) \omega^{123456} \\
0 = d(d\theta) \wedge \omega^3 \wedge \theta \wedge \omega^6 & = 4\!\left(ir' + p_4''\right)\omega^{123456}
\end{align*}
from which the result follows. $\lozenge$ \\

\indent A calculation using M\textsc{aple} shows that if the equations (5.7)-(5.8) of Lemma 5.2 hold, then the integrability conditions (5.1), (5.2), and (5.3) are all satisfied, and that
\begin{equation*}
d(dr) = \left( F\,\eta \wedge \overline{\eta} - 4ir''\phi \wedge \overline{\eta} \right) + dr'' \wedge \overline{\eta} \tag{5.9a}
\end{equation*}
where
\begin{equation*}
F = \textstyle -\left(  13 |p_4|^2 r + \frac{39}{2}r + 3r |r|^2 + 50i p_4^3 - \overline{p}_4r'' \right)\!. \tag{5.9b}
\end{equation*}
We summarize our discussion so far. \\

\noindent \textbf{Summary 5.3:} Nearly-K\"{a}hler structures of Type I are encoded by augmented coframings $((\eta, \overline{\eta}, \omega^3,$ $\theta, \overline{\theta}, \omega^6, \phi), (p,r), r'')$ on $P$ satisfying the following structure equations:  
\begin{align*}
d\eta & = i\phi \wedge \eta \tag{5.10a} \\
d\phi & = \textstyle i(6 + |r|^2 - |p|^2)\,\eta \wedge \overline{\eta} \\
\\
d\theta & = i\phi \wedge \theta + \left(ir\,\overline{\theta} \wedge \overline{\eta}  - \overline{p}\,\theta \wedge \eta \right) - 2p\,\text{Re}(\theta \wedge \overline{\eta}) - 4i\omega^6 \wedge \eta \\
d\omega^3 & = \textstyle \frac{5}{2}i\,\overline{\eta} \wedge \eta + \frac{3}{2}i\,\overline{\theta} \wedge \theta + 4\,\text{Re}(p\,\omega^3 \wedge \overline{\eta}) + 8\,\text{Re}(p\,\omega^6 \wedge \overline{\theta}) \\
d\omega^6 & = \textstyle \frac{3}{2}i\,\overline{\eta} \wedge \theta - \frac{3}{2}i\,\eta \wedge \overline{\theta} - 4\,\text{Re}(p\,\omega^6 \wedge \overline{\eta})
\end{align*}
and
\begin{align*}
dp & = \textstyle -\left( 6|p|^2 + \frac{3}{2} \right)\eta - \left(5p^2 -  i\overline{p}r \right) \overline{\eta} + ip\phi \tag{5.10b} \\
dr & = \textstyle -\left( 5ip^2 + \overline{p}r \right) \eta + r'' \overline{\eta} + 3ir\phi.
\end{align*}
where for ease of notation, we have set $p = p_4$. \\
\indent Augmented coframings satisfying the structure equations (5.10a)-(5.10b) will satisfy the integrability conditions (5.1), (5.2), (5.3), as well as (5.9a)-(5.9b). In the language of $\S$2.4, the functions $p$ and $r$ are the ``primary invariants," while $r''$ is the ``free derivative." \\

\noindent \textit{Remark:} The formulas of $\S$4.5.2 simplify considerably in the Type I setting.  In particular, we remark that the principal $G$-orbits in $M$ have mean curvature vector
\begin{equation*}
H = -4 \left(\text{Re}(p)\,e_1 + \text{Im}(p)\,e_2\right)
\end{equation*}
and have scalar curvature
$$\text{Scal} = \textstyle -\frac{9}{4} - 16\left|p\right|^2 < 0.$$
Thus, $p$ is essentially the mean curvature (or scalar curvature) of the principal orbits. $\Box$ \\

\noindent \textit{Remark:} From our remarks in $\S$4.5.3, we see that for Type I nearly-K\"{a}hler structures, the $2$-plane distribution normal to the principal $G$-orbits in $M$ is never integrable. $\Box$

\subsubsection{Local Existence/Generality}

\indent \indent We are now ready to state a local existence and generality theorem for Type I structures. \\

\noindent \textbf{Theorem 5.4:} Nearly-K\"{a}hler structures of Type I exist locally and depend on $2$ functions of $1$ variable in the sense of exterior differential systems.  In fact: \\
\indent For any $x \in \mathbb{R}^6$ and $(a_0, b_0) \in \mathbb{C}^2 \times \mathbb{C}$, there exists a Type I nearly-K\"{a}hler structure on an open neighborhood $U \subset \mathbb{R}^6$ of $x$ and an $\text{O}(2)$-coframe $f_x \in P|_x$ at $x$ for which
$$(p, r)(f_x) = a_0 \ \ \text{ and } \ \ r''(f_x) = b_0.$$ \\
\noindent \textit{Remark:} In a certain sense \cite{Bryant:2014aa}, the space of diffeomorphism classes of $k$-jets of Type I nearly-K\"{a}hler structures has dimension $2k + 4$. $\Box$ \\

\noindent \textit{Proof:} The discussion in $\S5.1.1$ shows that hypotheses (2.4) and (2.5) of Cartan's Third Theorem  (Theorem 2.2) are satisfied.  It remains to examine the tableau of free derivatives.  At a point $(u,v) \in \mathbb{R}^4 \times \mathbb{R}^2$, this is the vector subspace $A(u,v) \subset \text{Hom}(\mathbb{R}^7; \mathbb{R}^4) \cong \text{Mat}_{4 \times 7}(\mathbb{R})$ given by
$$A(u,v) = \left\{\left(\begin{array}{c c | c c c c c}
0 & 0 & 0 & 0 & 0 & 0 & 0 \\
0 & 0 & 0 & 0 & 0 & 0 & 0 \\ \hline
x & y & 0 & 0 & 0 & 0 & 0 \\
y & -x & 0 & 0 & 0 & 0 & 0 \end{array} \right) \colon x,y \in \mathbb{R} \right\}\!.$$
Since $A(u,v)$ is independent of the point $(u,v) \in \mathbb{R}^4 \times \mathbb{R}^2$, we can write $A = A(u,v)$ without ambiguity. We observe that $A$ is $2$-dimensional and has Cartan characters $\widetilde{s}_1 = 2$ and $\widetilde{s}_k = 0$ for $k \geq 2$.  One can check that $A$ is an involutive tableau, meaning that its prolongation $A^{(1)}$ satisfies $\dim(A^{(1)}) = 2 = \widetilde{s}_1 + 2\widetilde{s}_2 + \cdots + 7\widetilde{s}_7$. Thus, Cartan's Third Theorem applies, and we conclude the result. $\lozenge$ \\

\noindent \textit{Remark:} The complex characteristic variety of the tableau $A$ is
\begin{align*}
\Xi_A^{\mathbb{C}} & = \{[\xi] \in \mathbb{P}(\mathbb{C}^7) \colon w \otimes \xi \in A \text{ for some } w \in \mathbb{R}^4, w \neq 0\} \\
& = \{ [\xi] \in \mathbb{P}(\mathbb{C}^7) \colon (\xi_1 + i\xi_2)(\xi_1 - i\xi_2) = 0, \ \xi_3 = \cdots = \xi_7 = 0\}.
\end{align*}
\noindent The fact that the local generality of Type I structures is $2$ functions of $1$ variable, with complex characteristic variety consisting of two complex conjugate points, strongly suggests the possibility of a holomorphic interpretation of these structures.  $\Box$

\subsubsection{Incompleteness}

\indent \indent Nearly-K\"{a}hler structures of Type I cannot arise from a complete metric, as we now show.  Recall that the real Heisenberg group is the (non-compact) Lie group
$$\text{H}_3 = \left\{ \begin{pmatrix} 1 & x_1 & x_3 \\ 0 & 1 & x_2 \\ 0 & 0 & 1 \end{pmatrix} \colon x_i \in \mathbb{R} \right\} \leq \text{GL}_3(\mathbb{R}).$$ \\
\noindent \textbf{Proposition 5.5:} If $M$ is of Type I, then the universal cover of the acting Lie group $G$ is $\widetilde{G} = \text{H}_3 \times \mathbb{R}$.  In particular, the metric on $M$ is incomplete. \\

\noindent \textit{Proof:} As in the proof of Proposition 4.5, we identify $G$ with an integral $4$-fold of the  ideal $\mathcal{I}_G = \langle \omega^1, \omega^2, \phi \rangle = \langle \eta, \overline{\eta}, \phi \rangle$.  Under this identification, $\{\omega^3, \omega^4, \omega^5, \omega^6\}$ is a basis of left-invariant $1$-forms on $G$.  From (5.11), their exterior derivatives (mod $\mathcal{I}_G$) are given by
\begin{align*}
d\omega^3 & \equiv -3\omega^{45} - 8\,\text{Re}(p)\,\omega^{46} - 8\,\text{Im}(p)\,\omega^{65} \\
d\omega^4 & \equiv d\omega^5 \equiv d\omega^6 \equiv 0.
\end{align*}
\indent Let $\{X_3, X_4, X_5, X_6\}$ be a basis of $\mathfrak{g} = \{\text{left-invariant vector fields on }G\}$ whose dual basis is $\{\omega^3, \omega^4, \omega^5, \omega^6\}$. Let $Y = \frac{8}{3}\text{Im}(p)\,X_4 - \frac{8}{3}\text{Re}(p) X_5 + X_6$.  Then $\{X_3, X_4, X_5, Y\}$ is a basis of $\mathfrak{g}$ with
\begin{align*}
[X_4, X_5] & = 3X_3 & [X_3, X_4] & = 0 \\
[X_4, Y] & = 0 & [X_3, X_5] & = 0 \\
[X_5, Y] & = 0 & [X_3, Y] & = 0.
\end{align*}
This exhibits $\mathfrak{g}$ as the Lie algebra of the Lie group $\text{H}_3 \times \mathbb{R}$, and so the universal cover of $G$ is $\widetilde{G} = \text{H}_3 \times \mathbb{R}$.  Thus, Proposition 3.3 implies that the underlying metric is incomplete. $\lozenge$

\subsubsection{A Quasilinear Elliptic PDE System}

\indent \indent We now derive a quasilinear elliptic PDE system on the orbit space $\Sigma$ which Type I nearly-K\"{a}hler structures must satisfy.  For this, we make use of the setup of $\S$4.6, which we recall briefly.  Let $(z)$ be a local holomorphic coordinate on $\Sigma$, and write
$$\lambda \eta = \varpi^*(dz)$$
for some non-zero function $\lambda = e^{u+iv} \colon F \to \mathbb{C}$, where $u,v \colon F \to \mathbb{R}$.  Note that 
\begin{align*}
e^{2u}\,\eta \circ \overline{\eta} & = \varpi^*\!\left(dz \circ d\overline{z}\right).
\end{align*}
\indent In the Type I setting, we have that
\begin{align*}
\overline{p}\eta & = \varpi^*(\alpha) = \varpi^*(f\,dz) \\
\overline{r}\eta^3 & = \varpi^*(\beta) = \varpi^*(g\,dz^3)
\end{align*}
for some polynomial forms $\alpha \in \Gamma(K_\Sigma)$, $\beta \in \Gamma(K_\Sigma^{\otimes 3})$ and functions $f,g \colon \Sigma \to \mathbb{C}$. \\

\noindent \textbf{Proposition 5.6:} For any nearly-K\"{a}hler structure of Type I, the functions $f,g,u \colon \Sigma \to \mathbb{C}$ satisfy the following quasilinear elliptic PDE system:
\begin{align*}
\frac{\partial f}{\partial \overline{z}} & = -6|f|^2 - \frac{3}{2} e^{-2u} \tag{5.11} \\
\frac{\partial g}{\partial \overline{z}} & = 5if^2e^{-2u} - \overline{f}g \\
\Delta u & = -(6 + |g|^2 e^{6u} - |f|^2 e^{2u} )
\end{align*}
\noindent \textit{Proof:} We calculate
\begin{align*}
\overline{\partial}\alpha & = \overline{\partial}(f\,dz) = \nabla^{0,1}(f\,dz) = \frac{\partial f}{\partial \overline{z}}\,d\overline{z} \otimes dz \\
\varpi^*(\overline{\partial} \alpha) & = -\!\left(6p\overline{p} + \frac{3}{2} \right) \overline{\eta} \otimes \eta = \varpi^*\!\left( \left( -6|f|^2 - \frac{3}{2}e^{-2u} \right) d\overline{z} \otimes dz \right)
\end{align*}
and
\begin{align*}
\overline{\partial}\beta & = \overline{\partial}(g\,dz^3) = \nabla^{0,1}(g\,dz^3) = \frac{\partial g}{\partial \overline{z}}\,d\overline{z} \otimes dz^3 \\
\varpi^*(\overline{\partial} \beta) & = \left( 5i\overline{p}^2 - p\overline{r} \right) \overline{\eta} \otimes \eta^3 = \varpi^*\!\left( \left( 5i e^{-2u} f^2 - \overline{f}g \right) d\overline{z} \otimes dz^3 \right).
\end{align*}
This gives the first two equations.  For the last equation, we simply apply Lemma 4.7:
\begin{align*}
\Delta u = -K = -(6 + |r|^2 - |p|^2) = -(6 + e^{6u} \varpi^*|g|^2 - e^{2u} \varpi^*|f|^2).
\end{align*}
\noindent $\lozenge$ \\

\noindent \textit{Remark:} We have shown that a Type I structure determines a Riemann surface $\Sigma$ and functions $f,g,u \colon \Sigma \to \mathbb{C}$ which satisfy (5.11).  Conversely, given a Riemann surface $\Sigma$ and functions $f,g,u \colon \Sigma \to \mathbb{C}$ that solve (5.11), one can clearly recover $1$-forms $\eta, \overline{\eta}, \phi$ and functions $p, r$ satisfying
\begin{align*}
d\eta & = i\phi \wedge \eta \\
d\phi & = \textstyle i(6 + |r|^2 - |p|^2)\,\eta \wedge \overline{\eta} \\
dp & = \textstyle -\left( 6|p|^2 + \frac{3}{2} \right)\eta - \left(5p^2 -  i\overline{p}r \right) \overline{\eta} + ip\phi \\
dr & = \textstyle -\left( 5ip^2 + \overline{p}r \right) \eta + r'' \overline{\eta} + 3ir\phi.
\end{align*}
\indent We expect that, in fact, one can go further recover and $1$-forms $\omega^3, \theta, \overline{\theta}, \omega^6$ satisfying the remainder of the structure equations in (5.10a), and in this way recover the entire nearly-K\"{a}hler structure.  Showing this will require a deeper understanding of the structure equations (5.10a)-(5.10b). $\Box$

\subsection{Type II}

\indent \indent We now examine nearly-K\"{a}hler structures of Type II.  The integrability conditions for Type II structures are significantly more complicated than those for Type I.  To satisfy them, we will make a further frame adaptation and a change-of-variable. \\
\indent Ultimately, we will draw two conclusions.  First, we obtain (Theorem 5.10) a local existence/generality theorem for Type II structures.  Second, will show that the Lie group $G$ is solvable (Proposition 5.11), and hence that the underlying metrics are incomplete.

\subsubsection{A Frame Adaptation}

\indent \indent By definition, Type II structures are those with $\Vert q_1 \Vert = \Vert s_2 \Vert$ and $\Vert q_2 \Vert = \Vert s_1 \Vert$ and $\langle q_1, q_2 \rangle = \langle s_1, s_2 \rangle \neq 0$.  Thus, the $\text{O}(2)$-equivariant function $\frac{q_1}{s_2} = \frac{q_2}{s_1} \colon P \to \mathbb{C}$ maps into the unit circle $\mathbb{S}^1 \subset \mathbb{C}$.  Accordingly, we may adapt frames as follows: define the $\mathbb{Z}_2$-subbundle
$$P_1 = \{u \in P \colon q_1(u) = is_2(u)\} \subset P.$$
We refer to elements of $P_1$ as \textit{$\mathbb{Z}_2$-coframes}.  For the remainder of $\S$5.2, we work on $P_1$. \\
\indent The price we pay for this adaptation is the presence of additional torsion functions.  Indeed, on $P_1$ the $1$-form $\phi$ is no longer a connection form, but rather
$$\phi = \ell_1\omega^1 + \ell_2\omega^2$$
for some new $G$-invariant functions $\ell_1, \ell_2 \colon P_1 \to \mathbb{R}$.

\subsubsection{The Integrability Conditions}

\indent \indent We now move to solve the integrability conditions (5.1).  For this, we make the following change-of-variables.  Rather than work with $p_1, p_2, p_3, p_4, r \colon P_1 \to \mathbb{C}$, we will work with $t_1, \ldots, t_8$, $r_1, r_2 \colon P_1 \to \mathbb{R}$ defined by:
\begin{align*}
t_1 & = \text{Re}(p_1 + 4ip_4) & t_5 & = \text{Im}(p_1 + 4ip_4) & r_1 & = \text{Re}(r) \\
t_2 & = \textstyle \frac{1}{24}\,\text{Re}(p_3 + 4p_4) & t_6 & = \text{Im}(p_3 + 4p_4) & r_2 & = \text{Im}(r) \\
t_3 & = \text{Re}(p_2) & t_7 & = \text{Im}(p_2) & & \\
t_4 & = \text{Im}(p_4) & t_8 & = \text{Re}(p_4). & &
\end{align*}
The factor of $\frac{1}{24}$ appearing in $t_2$ is merely for the sake of clearing future denominators.  Since each $t_i$, $h_i$, and $\ell_i$ is $G$-invariant, we can write their exterior derivatives as
\begin{align*}
dt_i & = t_{i1} \omega^1 + t_{i2} \omega^2 & dh_i & = h_{i1}\omega^1 + h_{i2}\omega^2 & d\ell_i & = \ell_{i1}\omega^1 + \ell_{i2}\omega^2.
\end{align*}
We now state the Type II analogue of Lemma 5.2. \\

\noindent \textbf{Lemma 5.7:} Let $M$ be a nearly-K\"{a}hler manifold of Type II. \\
\indent (a) On the $\mathbb{Z}_2$-coframe bundle $P_1$, the following $12$ algebraic equations hold:
\begin{align*}
\text{Re}(q_1) & = -h_4 & t_5 & = 0 & h_2 & = -4t_2t_3 \tag{5.12} \\
\text{Im}(q_1) & = h_2 & t_6 & = t_1 + 8t_4 - 64 t_1 t_2^2 & h_3 & = -4t_1t_2 \\
\text{Re}(q_2) & = -h_3 & t_7 & = 0 & r_1 & = \ell_1 + t_4 \\
\text{Im}(q_2) & = h_1 + 3 & t_8 & = -t_2(2h_1 + 3) & r_2 & = \ell_2 + t_8 + 24t_2
\end{align*}
Thus, the torsion is expressible in terms of the $8$ real-valued functions
$$t_1, t_2, t_3, t_4 \ \ \ \ \text{ and } \ \ \ \ h_1, h_4, \ell_1, \ell_2.$$
\indent (b) The integrability conditions $d(d\omega^i) = 0$ are equivalent to the 12 algebraic equations (5.12) together with the equations
\small \begin{align*} 
t_{11} & = t_1(\ell_2 - 32h_1 t_2) & h_{11} & = H_{11}(t_1, t_2, t_3, t_4, h_1, h_4, \ell_1, \ell_2) \tag{5.13} \\
t_{12} & = t_1^2(64t_2^2 + 2) + t_1(4t_4 - \ell_1) + 2t_3^2 + 6(h_1 + 3) & h_{12} & = H_{12}(t_1, t_2, t_3, t_4, h_1, h_4, \ell_1, \ell_2) \\
t_{21} & = \textstyle 4t_2^2(2h_1 - 9)  - t_2\ell_2 - \frac{1}{2} & h_{41} & = H_{41}(t_1, t_2, t_3, t_4, h_1, h_4, \ell_1, \ell_2) \\
t_{22} & = t_2\ell_1 & h_{42} & = H_{42}(t_1, t_2, t_3, t_4, h_1, h_4, \ell_1, \ell_2)  \\
t_{31} & = 16t_2(h_4 t_1 - h_1 t_3) + \ell_2t_3 & \ell_{11} & = \textstyle u_1 + \frac{1}{2}G_1(t_1, t_2, t_3, t_4, h_1, h_4, \ell_1, \ell_2 ) \\
t_{32} & = t_3(192t_1 t_2^2 - 12t_4 - \ell_1) - 6h_4 & \ell_{12} & = \textstyle u_2 + \frac{1}{2}G_2(t_1, t_2, t_3, t_4, h_1, h_4, \ell_1, \ell_2 ) \\
t_{41} & = T_{41}(t_1, t_2, t_3, t_4, h_1, h_4, \ell_1, \ell_2) & \ell_{21} & = \textstyle u_2 - \frac{1}{2}G_2( t_1, t_2, t_3, t_4, h_1, h_4, \ell_1, \ell_2) \\
t_{42} & = T_{42}(t_1, t_2, t_3, t_4, h_1, h_4, \ell_1, \ell_2) & \ell_{22} & = \textstyle -u_1 + \frac{1}{2}G_1( t_1, t_2, t_3, t_4, h_1, h_4, \ell_1, \ell_2)
\end{align*} \normalsize
where we set $(u_1, u_2) = \left( \frac{1}{2}(\ell_{11} - \ell_{22}), \frac{1}{2}(\ell_{12} + \ell_{21})\right)$, and where the functions $T_{41}, T_{42}$ and $H_{11}, H_{12}$, $H_{41}, H_{42}$ and $G_1, G_2$ appearing on the right-hand sides of (5.13) are polynomial functions (of degree $\leq 5$) whose explicit formulas are listed in the appendix. \\

\noindent \textit{Proof:} (a) The left-most equations for $q_1$ and $q_2$ define our frame adaptation $P_1 \subset P$.  For the remaining eight equations, a calculation shows:
\begin{align*}
0 & = d(d\omega^4) \wedge \omega^{126} & \implies & & (h_1 + 3) t_7 + h_4 t_5  & = 0 \tag{5.14a} \\
0 & = d(d\omega^5) \wedge \omega^{126} & \implies & & h_2 t_5 - h_3 t_7 & = 0. \tag{5.14b}
\end{align*}
We rewrite (5.14a)-(5.14b) as
\begin{equation*}
\begin{pmatrix}
h_1 + 3 & h_4 \\
-h_3 & h_2 \end{pmatrix} \begin{pmatrix} t_7 \\ t_5 \end{pmatrix} = \begin{pmatrix} 0 \\ 0 \end{pmatrix}\!.
\end{equation*}
Since $M$ is of Type II, we have $(h_1 + 3)h_2 + h_3h_4 = \langle s_1, s_2 \rangle \neq 0$, from which it follows that $t_5 = t_7 = 0$.  Similarly, one can compute
\begin{align*}
0 &= d(d\omega^3) \wedge \omega^{126} & \implies & & 4 t_1 t_2 + h_3  & = 0 \tag{5.15a} \\
0 & = d(d\omega^3) \wedge \omega^{123} & \implies & & 4 t_2 t_3 + h_2  & = 0 \tag{5.15b} \\
0 & = d(d\omega^3) \wedge \omega^{236} & \implies & &  2 h_1 t_2 + 3 t_2 +  t_8 & = 0 \tag{5.16a} \\
0 & = d(d\omega^5) \wedge \omega^{123} & \implies & & 16 h_3 t_2 + t_1 + 8t_4 - t_6 & = 0 \tag{5.16b} \\
\\
0 & = d(d\omega^4) \wedge \omega^{235} & \implies & & h_4\!\left(\ell_1 - r_1 + t_4 \right) & = 0 \tag{5.17a} \\
0 & = d(d\omega^5) \wedge \omega^{235} & \implies & &  h_2(\ell_1 - r_1 + t_4)  & = 0 \tag{5.17b} \\
0 & = d(d\omega^4) \wedge \omega^{135} & \implies & &  h_4(\ell_2 - r_2 + 24t_2 + t_8)   & = 0 \tag{5.17c} \\
0 & = d(d\omega^5) \wedge \omega^{135} & \implies & &  h_2(\ell_2 - r_2 + 24t_2 + t_8)   & = 0. \tag{5.17d}
\end{align*}
Equations (5.15a)-(5.15b) now give the formulas for $h_2$ and $h_3$, while (5.16a)-(5.16b) give the formulas for $t_6$ and $t_8$.  Finally, since $M$ is of Type II, we have $(h_2)^2 + (h_4)^2 = \Vert s_2 \Vert^2 \neq 0$.  Thus, equations (5.17a)-(5.17d) give the remaining two equations.\\
\indent (b) This is a direct check of the equations remaining in $d(d\omega^i) = 0$.  $\lozenge$ \\

\indent A calculation using M\textsc{aple} shows that if the equations (5.12)-(5.13) of Lemma 5.7 hold, then $d(dt_i) = 0$ and $d(dh_i) = 0$ are also satisfied, and that
\begin{align*}
d(d\ell_1) & = F_1\,\omega^{12} + \left(du_1 \wedge \omega^1 + du_2 \wedge \omega^2\right) \tag{5.18} \\
d(d\ell_2) & = F_2\,\omega^{12} + \left(du_2 \wedge \omega^1 - du_1 \wedge \omega^2\right)\!,
\end{align*}
where $F_1$, $F_2$ are certain polynomial functions (of degree $\leq 8$) of $t_1, t_2, t_3, t_4$, $h_1, h_4$, $\ell_1, \ell_2$ and $u_1, u_2$ whose explicit formulas we will not list here. \\

\noindent \textbf{Summary 5.8:} Nearly-K\"{a}hler structures of Type II are encoded by augmented coframings $((\omega^1,$ $\ldots, \omega^6)$, $(t_1, t_2, t_3, t_4, h_1, h_4, \ell_1, \ell_2), (u_1, u_2))$ on the $\mathbb{Z}_2$-bundle $P_1 \to M$ satisfying the structure equations
\begin{align*}
d\omega^1 & = -\ell_1\,\omega^{12} \tag{5.19a} \\
d\omega^2 & = -\ell_2\,\omega^{12}
\end{align*}
and
\begin{align*}
d\omega^3 & = (2h_1 + 1)\,\omega^{12} - 4t_8 \omega^{13} + 2h_4 \omega^{15} - (t_1 + 4t_4)\,\omega^{23}  - 2h_2 \omega^{25} - t_3\omega^{26} \tag{5.19b} \\
& \ \ \ \ \ \ \ \ - t_3 \omega^{35} + 2h_2 \omega^{36}  - 3 \omega^{45}  - 24t_2 \omega^{46} - t_6 \omega^{56} \\
d\omega^4 & = (\ell_2 + 4t_8 + 24t_2)\,\omega^{14} - 2\ell_1 \omega^{15} + 2h_4\omega^{23} - \ell_1\omega^{24} - 2(\ell_2 + 12t_2)\, \omega^{25} + 2(h_1 + 1)\,\omega^{26} \\
& \ \ \ \ \ \ \ \ + 2(h_1 + 3)\,\omega^{35} - 8( t_8 - 3t_2)\,\omega^{36}   + 2h_4 \omega^{56} \\
d\omega^5 & = - (\ell_2 + 24t_2)\,\omega^{15} + 4 \omega^{16} - 2h_2 \omega^{23} - 24t_2 \omega^{24} + (\ell_1 + 4t_4)\,\omega^{25} + 2h_3 \omega^{26} \\
& \ \ \ \ \ \ \ \  + 2h_3 \omega^{35} + (t_6 - t_1 - 8t_4) \,\omega^{36}   - 2h_2 \omega^{56}  \\
d\omega^6 & = -2h_4 \omega^{12} + (2h_1 + 3)\,\omega^{15} - 4(t_8 - 6t_2)\,\omega^{16} - t_3\omega^{23} + 3\omega^{24} + 2h_3\,\omega^{25} + (t_6 - 4t_4)\,\omega^{26} \\
& \ \ \ \ \ \ \ \  + t_1\omega^{35}  - 2h_3\omega^{36} - t_3\omega^{56}
\end{align*} 
where $t_6, t_8, h_2, h_3$ are given by (5.12), and
\begin{align*}
dt_i & = t_{i1} \omega^1 + t_{i2} \omega^2 & dh_i & = h_{i1}\omega^1 + h_{i2}\omega^2 & d\ell_i & = \ell_{i1}\omega^1 + \ell_{i2}\omega^2 \tag{5.19c} 
\end{align*}
where $t_{11}, \ldots, t_{42}$ and $h_{11}, h_{12}, h_{41}, h_{42}$ and $\ell_{11}, \ell_{12}, \ell_{21}, \ell_{22}$ are given by (5.13). \\
\indent Augmented coframings satisfying the structure equations (5.19a)-(5.19c) and (5.12)-(5.13) will satisfy $d(d\omega^i) = 0$ and $d(dt_i) = d(dh_1) = d(dh_4) = 0$ and (5.18).   In the language of $\S$2.4, the functions $t_1, t_2, t_3, t_4$, $h_1, h_4$, $\ell_1, \ell_2$ are the ``primary invariants," while $u_1, u_2$ are the ``free derivatives."

\subsubsection{Local Existence/Generality}

\indent \indent We may now state the corresponding local existence/generality result for Type II structures. \\

\noindent \textbf{Theorem 5.9:} Nearly-K\"{a}hler structures of Type II exist locally and depend on $2$ functions of $1$ variable in the sense of exterior differential systems.  In fact: \\
\indent For any $x \in \mathbb{R}^6$ and any $(a_0, b_0) \in \mathbb{R}^8 \times \mathbb{R}^2$, there exists a Type II nearly-K\"{a}hler structure on an open neighborhood $U \subset \mathbb{R}^6$ of $x$ and a $
\mathbb{Z}_2$-coframe $f_x \in P_1|_x$ at $x$  for which
$$(t_1, t_2, t_3, t_4, h_1, h_4, \ell_1, \ell_2)(f_x) = a_0 \ \ \text{ and } \ \ (u_1, u_2)(f_x) = b_0$$ \\
\noindent \textit{Proof:} The discussion in $\S$5.2.2 shows that hypotheses (2.4) and (2.5) of Cartan's Third Theorem (Theorem 2.2) are satisfied.  It remains to examine the tableau of free derivatives.  At a point $(u,v) \in \mathbb{R}^8 \times \mathbb{R}^2$, this is the vector subspace $A(u,v) \subset \text{Hom}(\mathbb{R}^6; \mathbb{R}^8) \cong \text{Mat}_{8 \times 6}(\mathbb{R})$ given by
$$A(u,v) = \left\{\left(\begin{array}{c c | c c c c }
0 & 0 & 0 & 0 & 0 & 0  \\
\vdots & \vdots & \vdots & &  & \vdots  \\
0 & 0 & 0 & 0 & 0 & 0 \\ \hline
x & y & 0 & 0 & 0 & 0 \\
y & -x & 0 & 0 & 0 & 0  \end{array} \right) \colon x,y \in \mathbb{R} \right\}\!.$$
Since $A(u,v)$ is independent of the point $(u,v) \in \mathbb{R}^8 \times \mathbb{R}^2$, we can write $A = A(u,v)$ without ambiguity.  We observe that $A$ is $2$-dimensional and has Cartan characters $\widetilde{s}_1 = 2$ and $\widetilde{s}_k = 0$ for $k \geq 2$.  One can also check that $A$ is an involutive tableau, meaning that its prolongation $A^{(1)}$ satisfies $\dim(A^{(1)}) = 2 = \widetilde{s}_1 + 2\widetilde{s}_2 + \cdots + 6\widetilde{s}_6$. Thus, from Cartan's Third Theorem, we conclude the result. $\lozenge$

\subsubsection{Incompleteness}

\indent \indent As in the Type I setting, the non-compactness of the Lie group $G$ will prevent metrics of Type II from being complete. \\

\noindent \textbf{Proposition 5.10:} If $M$ is of Type II, then the Lie algebra $\mathfrak{g} = \text{Lie}(G)$ is solvable.  In particular, the metric on $M$ is incomplete. \\

\noindent \textit{Proof:} As in the proof of Proposition 4.5, we identify $G$ with an integral $4$-fold of the differential ideal $\mathcal{I}_G = \langle \omega^1, \omega^2, \phi \rangle = \langle \eta, \overline{\eta}, \phi \rangle$.  Under this identification, $\{\omega^3, \omega^4, \omega^5, \omega^6\}$ is a basis of $\mathfrak{g}^* = \{\text{left-invariant }1\text{-forms on }G\}$. \\
\indent Let $\zeta = \omega^5 + 8t_2\omega^6$, so that $\{\omega^3, \omega^4, \omega^6, \zeta\}$ is also a basis for $\mathfrak{g}^*$.  One can check that their exterior derivatives (mod $\langle \omega^1, \omega^2 \rangle$) are
\begin{align*}
d\omega^3 & \equiv -t_3\,\omega^3 \wedge \zeta - 3\,\omega^4 \wedge \zeta + t_6\,\omega^6 \wedge \zeta & d\omega^6 & \equiv t_1\,\omega^3 \wedge \zeta + t_3\,\omega^6 \wedge \zeta \\
d\omega^4 & \equiv 2(h_1 + 3)\,\omega^3 \wedge \zeta - 2h_4\,\omega^6 \wedge \zeta & d\zeta & \equiv 0,
\end{align*}
where we recall $t_6 = t_1 + 8t_4 - 64 t_1 t_2^2$. \\
\indent Let $\{X_3, X_4, X_5, Z\}$ be a basis of $\mathfrak{g} = \{\text{left-invariant vector fields on }G\}$ whose dual basis is $\{\omega^3, \omega^4, \omega^6, \zeta\}$.  Their Lie brackets are then
\begin{align*}
[X_3, Z] & = t_3X_3 - 2(h_1 + 3)X_4 - t_1 X_6 & [X_3, X_4] & = 0 \\
[X_4, Z] & = 3X_4 & [X_3, X_6] & = 0 \\
[X_6, Z] & = -t_6 X_3 + 2h_4 X_4 - t_3X_6 & [X_4, X_6] & = 0.
\end{align*}
From this, it is clear that $[[\mathfrak{g}, \mathfrak{g}], [\mathfrak{g}, \mathfrak{g}]] = 0$, so that $\mathfrak{g}$ is solvable. (Note, however, that $\mathfrak{g}$ is not nilpotent in general.) Thus, Proposition 3.3 implies that the underlying metric is incomplete. $\lozenge$

\subsection{Type III}

\indent \indent We now consider nearly-K\"{a}hler structures of Type III.  This is perhaps the most interesting case, as there is the possibility for complete metrics to exist in this class.  Unfortunately, the integrability conditions (5.1)-(5.4) are even more complicated than those of Type II. \\
\indent Examining these conditions leads us to several changes-of-variable ($\S$5.3.1, $\S$5.3.2).  The upshot is that the intrinsic torsion
$$(p_1, p_2, p_3, p_4, q_1, q_2, r; h_1, h_2, h_3, h_4) \colon P \to \mathbb{C}^7 \oplus \mathbb{R}^4$$
will be recast as a function
$$(u, v, z, r;\, t_0, t_1, t_2) \colon P \to \mathbb{C}^4 \oplus \mathbb{R}^3.$$
Even with this repackaging, however, we find it difficult to solve (5.1) in general.  As such, we will impose the ansatz that $\mathfrak{g} = \mathfrak{su}(2) \oplus \mathfrak{u}(1)$ (which by Proposition 3.3 is the case of most interest anyway).  This will let us normalize the function $t_2$, which simplifies (5.1) further. \\
\indent Finally, after this normalization, we are able to solve all the integrability conditions by means of a computer algebra system (Lemma 5.14), thus yielding a local existence/generality result (Theorem 5.15) for the nearly-K\"{a}hler structures with $\mathfrak{g} = \mathfrak{su}(2) \oplus \mathfrak{u}(1)$. 

\subsubsection{A First Change-of-Variable}

\indent \indent By definition, Type III structures are those with $\langle q_1, q_2 \rangle = \langle s_2, s_1 \rangle = 0$ and $q_1$, $q_2$, $s_1$, $s_2$ not all zero.  Recalling also that $q_1s_1 - q_2s_2 = 0$, we have:
\begin{align*}
& \text{rank}\!\begin{pmatrix} q_1 & q_2 \\ s_2 & s_1 \end{pmatrix} = 1, \text{ and } \\
& q_1\overline{q}_2 \text{ and } s_2\overline{s}_1 \text{ both pure imaginary}.
\end{align*}
This leads us to factor
$$\begin{pmatrix} q_1 & q_2 \\ s_2 & s_1 \end{pmatrix} = \begin{pmatrix} z \\ w \end{pmatrix} \begin{pmatrix} it_1 & t_2 \end{pmatrix} =
\begin{pmatrix} it_1 z & t_2 z \\ it_1 w & t_2w \end{pmatrix}\!,$$
where $z,w \colon P \to \mathbb{C}$ and $t_1, t_2 \colon P \to \mathbb{R}$.  Note that by definition of Type III, we cannot have $(z,w) = (0,0)$, nor can we have $(t_1, t_2) = (0,0)$. \\

\noindent \textit{Remark:} We caution that the functions $z, w, t_1, t_2$ are not uniquely defined: we may replace $(z, w, it_1, t_2)$ with $(cz, cw, it_1/c, t_2/c)$ for any non-vanishing function $c: P \to \mathbb{R}$. $\Box$

\subsubsection{A Second Change-of-Variable}

\indent \indent We now solve $4[\binom{4}{3} + \binom{3}{3} + \binom{3}{3}] = 24$ of the $4 \binom{6}{3} = 80$ equations arising in (5.1).  Namely, we will solve $d(d\nu) \wedge \eta \wedge \overline{\eta} = 0$ and $d(d\nu) \wedge \eta \wedge \omega^{36} = 0$ and $d(d\nu) \wedge \overline{\eta} \wedge \omega^{36} = 0$ for $\nu \in \{\omega^3, \theta, \overline{\theta}, \omega^6\}$.  This is accomplished in the following lemma: \\  

\noindent \textbf{Lemma 5.11:} Let $M$ be a nearly-K\"{a}hler manifold of Type III. \\
\indent (a) There exist functions $u, v, \widehat{v} \colon P \to \mathbb{C}$ such that:
\begin{align*}
p_1 + 4ip_4 & = 6t_2u & p_3 + 4p_4 & = 24v. \\
p_2 & = -6t_1u & -p_3 + 4p_4 & = 24\widehat{v}.
\end{align*}
\indent (b) On the $\text{O}(2)$-coframe bundle $P$, the following algebraic relations hold:
\begin{align*}
\text{Im}(w) & = -24\,\text{Re}(u\overline{v}) \tag{5.20} \\
12\widehat{v} & = z(t_1^2\overline{u} + 4it_2\overline{v}) - iw(t_1^2u - 4it_2v) - 3it_2u \tag{5.21}
\end{align*} \\
\noindent \textit{Proof:} The existence of $v, \widehat{v}$ is immediate.  Let us set $y = p_1 + 4ip_4$ and expand the identities $d(d\omega^3) \wedge \eta \wedge \overline{\eta} = 0$ and $d(d\omega^6) \wedge \eta \wedge \overline{\eta} = 0$.  This yields, for example,
\begin{align*}
d(d\omega^6) \wedge \eta \wedge \overline{\eta} \wedge \theta & = 0 & \implies &  & is_2y - q_1\overline{y} - s_1p_2 - iq_2 \overline{p}_2 & = 0 \tag{5.22a} \\
d(d\omega^6) \wedge \eta \wedge \overline{\eta} \wedge \overline{\theta} & = 0 &   \implies & & \overline{q}_1y + i \overline{s}_2\overline{y} - i\overline{q}_2p_2 + \overline{s}_1 \overline{p}_2 & = 0 \tag{5.22b}
\end{align*}
and
\begin{align*}
d(d\omega^3) \wedge \eta \wedge \overline{\eta} \wedge \omega^3 & = 0 & \implies &  & 4(\overline{p}_2v + p_2\overline{v}) & = -(s_2 + \overline{s}_2) \tag{5.23a} \\
d(d\omega^6) \wedge \eta \wedge \overline{\eta} \wedge \omega^3 & = 0 & \implies &  & 4(\overline{y}v + y\overline{v})  & = i(s_1 - \overline{s}_1) \tag{5.23b} \\
d(d\omega^6) \wedge \eta \wedge \overline{\eta} \wedge \omega^6 & = 0 & \implies &  & \overline{p}_2y - p_2\overline{y} & = 0. \tag{5.24}
\end{align*}
In light of (5.24), we see that in order for the linear system (5.23a)-(5.23b) to have solutions, we must have
\begin{equation*}
i(s_1 - \overline{s}_1) \begin{pmatrix} \overline{p}_2 \\ p_2 \end{pmatrix} + (s_2 + \overline{s}_2) \begin{pmatrix} \overline{y} \\ y \end{pmatrix} = \begin{pmatrix} 0 \\ 0 \end{pmatrix}\!. \tag{5.25}
\end{equation*}
Regard the equations (5.22a), (5.22b), and (5.25) as a homogeneous linear system:
\begin{equation*}
\begin{pmatrix}
is_2 & -q_1 & -s_1 & -iq_2 \\
\overline{q}_1 & i\overline{s}_2 & -i\overline{q}_2 & \overline{s}_1 \\
s_2 + \overline{s}_2 & 0 & i(s_1 - \overline{s}_1) & 0 \\
0 & s_2 + \overline{s}_2 & 0 & i(s_1 - \overline{s}_1)
\end{pmatrix}
\begin{pmatrix}
y \\
\overline{y} \\
p_2 \\
\overline{p}_2
\end{pmatrix} = \begin{pmatrix} 0 \\ 0 \\ 0 \\ 0 \end{pmatrix}
\end{equation*}
The solutions to this system are of the form
\begin{equation*}
\begin{pmatrix}
y \\ \overline{y} \\ p_2 \\ \overline{p}_2 \end{pmatrix} =
6u \begin{pmatrix} t_2 \\ 0 \\ -t_1 \\ 0 \end{pmatrix} + 6 \overline{u} \begin{pmatrix} 0 \\ t_2 \\ 0 \\ -t_1 \end{pmatrix}
\end{equation*}
for some $u \colon P \to \mathbb{C}$.  This proves (a).  The only equations left in $d(d\nu) \wedge \eta \wedge \overline{\eta} = 0$, $d(d\nu) \wedge \eta \wedge \omega^{36} = 0$ and $d(d\nu) \wedge \overline{\eta} \wedge \omega^{36} = 0$ for $\nu \in \{\omega^3, \theta, \overline{\theta}, \omega^6\}$ are exactly those in the statement of (b). $\lozenge$ \\

\noindent \textbf{Bookkeeping 5.12:} We pause to unwind the notational changes.  By definition of Type III and by Lemma 5.11(a), our torsion functions $(p_i, q_i, r, s_i)$ are now expressed as follows:
\begin{align*}
p_1 & = 6t_2u - 12i(v + \widehat{v}) & p_3 & = 12(v - \widehat{v}) & r & = r & q_1 & = it_1z & s_1 & = t_2w &  \\
p_2 & = -6t_1 u & p_4 & = 3(v + \widehat{v}) & &  & q_2 & = t_2z & s_2 & = it_1w.
\end{align*}
That is, the torsion is expressed in terms of
$$u,v,\widehat{v}, z, w, r \colon P \to \mathbb{C} \ \ \ \text{ and } \ \ \ t_1, t_2 \colon P \to \mathbb{R}.$$
\indent By Lemma 5.11(b), the functions $\widehat{v}$ and $\text{Im}(w)$ can be expressed in terms of the others.  Hence, setting $t_0 = \text{Re}(w)$, we regard the torsion as a function
$$(u, v, z, r;\, t_0, t_1, t_2) \colon P \to \mathbb{C}^4 \oplus \mathbb{R}^3.$$
One can check that these functions are $\text{O}(2)$-equivariant.  Indeed, modulo $\langle \omega^1, \omega^2 \rangle = \langle \eta, \overline{\eta}\rangle$:
\begin{align*}
du & \equiv iu\phi & dz & \equiv 2iz\phi & dt_1 & \equiv 0 \tag{5.26} \\
dv & \equiv iv\phi & dr & \equiv 3ir\phi & dt_2 & \equiv 0 \\
 &  &  &  & dt_0 & \equiv 0.
\end{align*}
\subsubsection{The Ansatz $\mathfrak{g} = \mathfrak{su}(2) \oplus \mathfrak{u}(1)$}

\indent \indent Having solved the $24$ of the $80$ equations arising in (5.1), we aim to solve the remaining $56$ of them.  To this end, we restrict attention to the case where $G$ is a finite quotient of $\text{SU}(2) \times \text{U}(1)$.  This ansatz imposes inequalities on the torsion which allow us to normalize $t_2$. \\

\noindent \textbf{Lemma 5.13:} Let $M$ be of Type III.  If $\mathfrak{g} = \mathfrak{su}(2) \oplus \mathfrak{u}(1)$, then the (real) $1$-form
$$\sigma_6 := 3\left(iu\overline{z} - \overline{u}\overline{w}\right) \theta - 3\left(i\overline{u}z + uw\right) \overline{\theta} + \left( |w|^2 - |z|^2 \right)\omega^6$$
is non-vanishing, and the (real) symmetric matrix
\begin{equation*}
Q := \begin{pmatrix}
\frac{1}{3}t_2 & -2t_1\,\text{Im}(u) & -2t_1\,\text{Re}(u) \\
-2t_1\,\text{Im}(u) & 48\,\text{Im}(u)\,\text{Re}(v) - \text{Im}(z) + t_0 & 24\,\text{Re}(uv) - \text{Re}(z) \\
-2t_1\,\text{Re}(u) & 24\,\text{Re}(uv) - \text{Re}(z) & -48\,\text{Re}(u)\,\text{Im}(v) + \text{Im}(z) + t_0
\end{pmatrix}
\end{equation*}
is positive-definite or negative-definite.  In particular, $t_2$ is nowhere-vanishing.  \\ 

\noindent \textit{Proof:} As in the proof of Proposition 4.5, we identify $G$ with an integral $4$-fold of the differential ideal $\mathcal{I}_G = \langle \omega^1, \omega^2, \phi \rangle = \langle \eta, \overline{\eta}, \phi \rangle$.  Under this identification, $\{\omega^3, \omega^4, \omega^5, \omega^6\}$ is a basis of $\mathfrak{g}^* = \{\text{left-invariant }1\text{-forms on }G\}$. \\
\indent Suppose that $\mathfrak{g} = \mathfrak{su}(2) \oplus \mathfrak{u}(1)$.  Then $\mathfrak{g}$ has a non-zero center, so there exists a non-zero element of $\mathfrak{g}^*$ which is closed.  A calculation shows that the only elements of $\mathfrak{g}^*$ which are closed are multiples of
$$\sigma_6 := 3\left(iu\overline{z} - \overline{u}\overline{w}\right) \theta - 3\left(i\overline{u}z + uw\right) \overline{\theta} + \left( |w|^2 - |z|^2 \right)\omega^6$$
Thus, $\sigma_6$ is non-vanishing. \\
\indent We now observe that $\{\sigma_3, \sigma_4, \sigma_5, \sigma_6\}$ is a basis for $\mathfrak{g}^*$, where we are defining
\begin{align*}
(\sigma_3, \sigma_4, \sigma_5) & := (t_2\omega^3 - t_1\omega^6, \, \text{Im}(\chi), \, \text{Re}(\chi)) & \chi & := -3(2t_1u\,\omega^3 - \theta - 8iv\,\omega^6).
\end{align*}
One can calculate that modulo $\mathcal{I}_G$,
\begin{align*}
d\!\begin{pmatrix} \sigma_3 \\ \sigma_4 \\ \sigma_5 \end{pmatrix} & \equiv Q
\begin{pmatrix}
\sigma_4 \wedge \sigma_5 \\
\sigma_5 \wedge \sigma_3 \\
\sigma_3 \wedge \sigma_4
\end{pmatrix} \\
d\sigma_6 & \equiv 0.
\end{align*}
Since $\{\sigma_3, \sigma_4, \sigma_5\}$ is a basis of $\mathfrak{su}(2)^*$, this coefficient matrix $Q$ must be positive-definite or negative-definite, and hence $t_2$ is nowhere-vanishing. $\lozenge$

\subsubsection{Local Existence/Generality}

\indent \indent We now move to solve the integrability conditions (5.1) in the case of $\mathfrak{g} = \mathfrak{su}(2) \oplus \mathfrak{u}(1)$.  By Lemma 5.13, the function $t_2$ is nowhere-vanishing.  Recalling that $z,w, t_1, t_2$ are only defined up to scaling by a nowhere-vanishing function $c \colon P \to \mathbb{R}$, we shall choose $c$ so that
$$t_2 = 1.$$
Thus, the torsion of the $\text{O}(2)$-structure is now encoded by $(u, v, z, r; t_0, t_1) \colon P \to \mathbb{C}^4 \oplus \mathbb{R}^2$.
Since $u, v, z, r$, $t_0, t_1$ are $G$-invariant and $\text{O}(2)$-equivariant (by (5.26)), we may express their exterior derivatives as
\begin{align*}
du & = u'\eta + u''\overline{\eta} + iu\phi & dz & = z'\eta + z''\overline{\eta} + 2iz\phi & dt_1 & = t_1'\eta + t_1''\overline{\eta} \\
dv & = v'\eta + v''\overline{\eta} + iv\phi & dr & = r'\eta + r''\overline{\eta} + 3ir\phi  & dt_0 & = t_0'\eta + t_0''\overline{\eta}.
\end{align*}
\noindent The Type III analogue of Lemma 5.2(b) and Lemma 5.7(b) is the following: \\

\noindent \textbf{Lemma 5.14:} Let $M$ be a Type III nearly-K\"{a}hler structure with $\mathfrak{g} = \mathfrak{su}(2) \oplus \mathfrak{u}(1)$.  With the equations of Lemma 5.11 imposed and with the normalization $t_2 = 1$ in place, the integrability conditions (5.1) are equivalent to
\begin{align*}
u' & = f_1(u, \overline{u}, v, \overline{v}, z, \overline{z}, r, \overline{r}, t_0, t_1) & u'' & = f_2(u, \overline{u}, v, \overline{v}, z, \overline{z}, r, \overline{r}, t_0, t_1) \tag{5.27} \\
v' & = f_3(u, \overline{u}, v, \overline{v}, z, \overline{z}, r, \overline{r}, t_0, t_1) & v'' & = f_4(u, \overline{u}, v, \overline{v}, z, \overline{z}, r, \overline{r}, t_0, t_1) \\
z' & = f_5(u, \overline{u}, v, \overline{v}, z, \overline{z}, r, \overline{r}, t_0, t_1) & z'' & = f_6(u, \overline{u}, v, \overline{v}, z, \overline{z}, r, \overline{r}, t_0, t_1) \\
r' & = f_7(u, \overline{u}, v, \overline{v}, z, \overline{z}, r, \overline{r}, t_0, t_1) &  & \\
t_0' & = f_8(u, \overline{u}, v, \overline{v}, z, \overline{z}, r, \overline{r}, t_0, t_1) & t_0'' & = f_{9}(u, \overline{u}, v, \overline{v}, z, \overline{z}, r, \overline{r}, t_0, t_1) \\
t_1' & = f_{10}(u, \overline{u}, v, \overline{v}, z, \overline{z}, r, \overline{r}, t_0, t_1) & t_1'' & = f_{11}(u, \overline{u}, v, \overline{v}, z, \overline{z}, r, \overline{r}, t_0, t_1) \\
K & = f_{12}(u, \overline{u}, v, \overline{v}, z, \overline{z}, r, \overline{r}, t_0, t_1). & &
\end{align*}
for certain functions $f_1, \ldots, f_{12}$ of $u, \overline{u}, v, \overline{v}, z, \overline{z}, r, \overline{r}, t_0, t_1$. \\

\noindent \textit{Remark:} The explicit expressions for $f_1, \ldots, f_{12}$ are sufficiently cumbersome that we will not list them here.  They turn out to be polynomial functions of degree $\leq 10$. $\Box$ \\

\noindent \textit{Proof:} With the equations of Lemma 5.11 imposed, and with $t_2 = 1$ imposed, there are $56 = 80 - 24$ polynomial equations remaining in (5.1) involving the functions
$$u, \overline{u}, v, \overline{v}, z, \overline{z}, r, \overline{r}, t_0, t_1 \ \ \text{ and } \ \ u', u'', \overline{u}', \overline{u}'', v', v'', \overline{v}', \overline{v}'', z', z'', \overline{z}', \overline{z}'', r', r'', \overline{r}', \overline{r}'', t_0', t_0'', t_1', t_1'', K.$$
A direct application of a computer algebra system (we used M\textsc{aple}) will solve these 56 equations, yielding lengthy explicit formulas for $u, u'', \ldots, t_1', t_1'', K$ in terms of $u, \overline{u}, v, \overline{v}, z, \overline{z}, r, \overline{r}, t_0, t_1$.  $\lozenge$ \\

\indent As in the Type I and Type II settings, a calculation using M\textsc{aple}  now shows that if the equations (5.27) of Lemma 5.14 hold, then $d(d\phi) = 0$ and $d(du) = d(dv) = d(dz) = d(dt_0) = d(dt_1) = 0$ are all satisfied, and that
\begin{align*}
d(dr) & = \left( F\,\eta \wedge \overline{\eta} - 4ir''\phi \wedge \overline{\eta}\right) +  dr'' \wedge \overline{\eta}
\end{align*}
where $F$ is a certain polynomial function (of degree $14$) of $u, \overline{u}, v, \overline{v}, z, \overline{z}, r, \overline{r}, t_0, t_1$ whose explicit formula we will not list here. \\

\indent The upshot of this discussion is that the integrability conditions (2.4) and (2.5) of Cartan's Third Theorem are finally satisfied.  In particular, we obtain the following local existence/generality result: \\

\noindent \textbf{Theorem 5.15:} Nearly-K\"{a}hler structures (of Type III) for which $G$ is a finite quotient of $\text{SU}(2) \times \text{U}(1)$ exist locally and depend on $2$ functions of $1$ variable in the sense of exterior differential systems.  In fact: \\
\indent For any $x \in \mathbb{R}^6$ and $(a_1, a_2, a_3, a_4, \widetilde{a}_5, \widetilde{a}_6, b_0) \in \mathbb{C}^4 \times \mathbb{R}^2 \times \mathbb{C}$ with $(a_3, \widetilde{a}_5, \text{Re}(a_1 \overline{a}_2)) \neq (0,0,0)$, there exists a (Type III) nearly-K\"{a}hler structure with $\mathfrak{g} = \mathfrak{su}(2) \oplus \mathfrak{u}(1)$ on an open neighborhood $U \subset \mathbb{R}^6$ of $x$ and an $\text{O}(2)$-coframe $f_x \in P|_x$ at $x$ for which
$$(u,v,z,r)(f_x) = (a_1, a_2, a_3, a_4) \ \ \text{ and } \ \ (t_0, t_1)(f_x) = (\widetilde{a}_5, \widetilde{a}_6) \ \ \text{ and } \ \ r''(f_x) = b_0.$$ \\
\noindent \textit{Remark:} The unusual looking requirement $(a_3, \widetilde{a}_5, \text{Re}(a_1 \overline{a}_2)) \neq (0,0,0)$ is simply the condition $(z, \text{Re}(w), \text{Im}(w)) \neq (0,0, 0)$ mentioned at the start of $\S$5.3.1.  That is, it is exactly the condition $(q_1, q_2, s_1, s_2) \neq (0,0, 0, 0)$ forming part of the definition of ``Type III." $\Box$ \\

\noindent \textit{Proof:} It remains only to examine the tableau of free derivatives.  This proceeds exactly as in the cases of Types I and II, so we omit the details. $\lozenge$

\section{A One-Parameter Family of Incomplete Cohomogeneity-One Metrics}

\indent \indent With the issue of local existence settled, it is natural to seek special solutions.  That is, one would like to find explicit examples of nearly-K\"{a}hler structures with $\mathfrak{g} = \mathfrak{su}(2) \oplus \mathfrak{u}(1)$ with one or more of the invariants $(u, v, z, r, t_0, t_1)$ prescribed in advance.  For example, it is natural to ask that the intrinsic torsion
$$(u, v, z, r; t_0, t_1) \colon P \to \mathbb{C}^4 \oplus \mathbb{R}^2$$
take values in a prescribed $\text{O}(2)$-invariant subset of $\mathbb{C}^4 \oplus \mathbb{R}^2$. \\
\indent In this section, we comment on those nearly-K\"{a}hler structures with $\mathfrak{g} = \mathfrak{su}(2) \oplus \mathfrak{u}(1)$ that satisfy $u = 0$.  A priori, it is not at all clear how many such structures exist (or that any exist at all). \\

\indent Let $M$ be a nearly-K\"{a}hler $6$-manifold with $\mathfrak{g} = \mathfrak{su}(2) \oplus \mathfrak{u}(1)$ satisfying $u = 0$.  The condition $u = 0$ turns out to be quite strong: a calculation using M\textsc{aple} shows that the integrability conditions resulting from $u = 0$ imply that $z = 0$, $r = 0$, and $t_0 = 3$.  Thus, the torsion is encoded in the function
$$A = (a_1, a_2, a_0) := (24\,\text{Re}(v), 24\,\text{Im}(v), 3t_1) \colon P \to \mathbb{R}^3.$$
\noindent \textit{Remark:} Unwinding the various changes-of-variable, we find that we are exactly in the situation of
\begin{align*}
(p_1, p_2, p_3, p_4) & = (0, 0, 24v, 0) & r & = 0 \\
(q_1, q_2) & = (0,0) & (h_1, h_2, h_3, h_4) & = (0, 0, 0, 3t_1).
\end{align*}
That is, $(a_1, a_2, a_0) = (\text{Re}(p_3), \text{Im}(p_3), h_4)$. $\Box$ \\

\indent The structure equations on $P$ now read as follows:
\begin{align*}
d\omega^1 & = -\phi \wedge \omega^2 & d\omega^3 & = \omega^{12} + a_0 \omega^{15} - a_0 \omega^{24}  - 3\omega^{45} - a_1\omega^{46} - a_2 \omega^{56} \tag{6.1} \\
d\omega^2 & = \phi \wedge \omega^1 & d\omega^4 & = a_0\omega^{23} - \omega^{26} + 3\omega^{35}  + a_1\omega^{36} + a_0\omega^{56} -\phi \wedge \omega^5 \\
d\phi & = 3( a_0^2 + 1)\, \omega^{12} & d\omega^5 & = - a_0 \omega^{13} + \omega^{16} - 3\omega^{34} + a_2\omega^{36} - a_0\,\omega^{46} + \phi \wedge \omega^4 \\
& &  d\omega^6 & = -2a_0\,\omega^{12} + a_1 \omega^{16} + a_2 \omega^{26}
\end{align*}
and
\begin{align*}
dA = \begin{pmatrix} da_1 \\ da_2 \\ da_0 \end{pmatrix} =
\begin{pmatrix} 
a_0^2 - a_1^2 - 3 & -a_1 a_2 & -a_2 \\
-a_1a_2 & a_0^2 - a_2^2 - 3 & a_1 \\
-2a_0 a_1 & -2a_0 a_2 & 0 \end{pmatrix}
\begin{pmatrix} \omega^1 \\ \omega^2 \\ \phi \end{pmatrix}\!.
\tag{6.2}
\end{align*}
These equations satisfy $d(d\omega^i) = d(d\phi) = 0$ and $d(da_i) = 0$.  Thus, Cartan's Third Theorem yields: \\

\noindent \textbf{Theorem 6.1:} Nearly-K\"{a}hler structures of Type III for which $\mathfrak{g} = \mathfrak{su}(2) \oplus \mathfrak{u}(1)$ and $u = 0$ exist locally.  In fact: \\
\indent For any $(x_0, y_0, z_0) \in \mathbb{R}^3$, there exists a Type III nearly-K\"{a}hler structure with $\mathfrak{g} = \mathfrak{su}(2) \oplus \mathfrak{u}(1)$ and $u = 0$ on an open neighborhood $U \subset \mathbb{R}^6$ of $x$ and an $\text{O}(2)$-coframe $f_x \in P|_x$ at $x$ for which
$$(a_1, a_2, a_0)(f_x) = (x_0, y_0, z_0).$$ \\
\indent We can completely describe the induced metric on the principal orbits: \\

\noindent \textbf{Proposition 6.2:} Let $M$ be a Type III nearly-K\"{a}hler $6$-manifold with $\mathfrak{g} = \mathfrak{su}(2) \oplus \mathfrak{u}(1)$ and $u = 0$.  Then the principal $G$-orbits in $M$ are locally isometric to $\mathbb{S}^3(\frac{2}{3}) \times \mathbb{S}^1$. \\

\noindent \textit{Proof:} Let $N^4$ be a principal $G$-orbit.  We will work on the $\text{O}(2)$-bundle given by $\pi^{-1}(N) \to N$.  Note that on $\pi^{-1}(N)$, the functions $a_0, a_1, a_2$ are all constant, and that $\omega^1 = \omega^2 = 0$.  Therefore, on $\pi^{-1}(N)$, the structure equations (6.1) read:
\begin{align*}
d\begin{pmatrix} \omega^3 \\ \omega^4 \\ \omega^5 \\ \omega^6 \end{pmatrix} & = -\begin{pmatrix}
0 & \frac{3}{2}\omega^5 + a_1 \omega_6 & -\frac{3}{2}\omega^4 + a_2 \omega^6 & 0 \\
-\frac{3}{2}\omega^5 - a_1\omega^6 & 0 & \frac{3}{2}\omega^3 - a_0 \omega^6 - \phi & 0 \\
\frac{3}{2}\omega^4 - a_2\omega^6 & -\frac{3}{2}\omega^3 + a_0 \omega^6 + \phi & 0 & 0 \\
0 & 0 & 0 & 0
\end{pmatrix} \wedge \begin{pmatrix} \omega^3 \\ \omega^4 \\ \omega^5 \\ \omega^6 \end{pmatrix} \tag{6.3} \\
d\phi & = 0.
\end{align*}
Since $d\omega^6 = 0$, locally it is the case that $\omega^6 = df$ for some submersion $f \colon N^4 \to \mathbb{S}^1$.  The equations (6.3) show that the $f$-horizontal distribution $\langle \omega^3, \omega^4, \omega^5 \rangle$ is integrable.  Moreover, a calculation using (6.3) shows that the $f$-fibers are totally-geodesic and have constant curvature $\frac{9}{4}$.  Thus, $N$ is locally isometric to a product $\mathbb{S}^3(\frac{2}{3}) \times \mathbb{S}^1$. $\lozenge$ \\

\indent Before stating our main result on Type III nearly-K\"{a}hler structures with $u = 0$, we make two preliminary observations. \\
\indent First, note that if $A \colon P \to \mathbb{R}^3$ has rank $r$, then $\dim(\text{Aut}_{\text{O}(2)}) = 7 - r$.  Indeed, on the one hand, $\text{Aut}_{\text{O}(2)}$ acts freely on $P$.  On the other hand, it is a general fact \cite{MR1062197} that $\text{Aut}_{\text{O}(2)}$ acts transitively on the level sets of $A$, which by hypothesis are $(7-r)$-dimensional. \\
\indent Second, note that the substitution
\begin{align*}
\beta_3 & := 3 \omega^3 - a_0 \omega^6 - \phi \tag{6.4} \\
\beta_4 & := 3\omega^4 - a_2 \omega^6 - a_0 \omega^1 \\
\beta_5 & := 3\omega^5 + a_1 \omega^6 - a_0 \omega^2
\end{align*}
greatly simplifies the structure equations (6.1):
\begin{align*}
d\omega^1 & = - \phi \wedge \omega^2 & d\beta_3 & = -\beta_4 \wedge \beta_5 \tag{6.5} \\
d\omega^2 & = \phi \wedge \omega^1 & d\beta_4 & = -\beta_5 \wedge \beta_3 \\
d\phi & = 3(a_0^2 + 1)\,\omega^{12} & d\beta_5 & = -\beta_3 \wedge \beta_4 \\
d\omega^6 & = -2a_0\omega^{12} + (a_1 \omega^1 + a_2\omega^2) \wedge \omega^6.
\end{align*}
Our result is: \\

\noindent \textbf{Theorem 6.3:} Let $M$ be a Type III nearly-K\"{a}hler $6$-manifold with $\mathfrak{g} = \mathfrak{su}(2) \oplus \mathfrak{u}(1)$ and $u = 0$. \\
\indent (a) The image of $A \colon P \to \mathbb{R}^3$ is one of the following: \\
\indent \indent (i) One of the points $(0, 0, \pm \sqrt{3})$. \\
\indent \indent (ii) An open subset of the $2$-sphere $\{a_1^2 + a_2^2 + (a_0 - c)^2 = c^2 - 3\}$, where $\sqrt{3} < |c| < \infty$. \\
\indent \indent (iii) An open subset of the $2$-plane $\{a_0 = 0\}$. \\
\indent (b) In case (i), the nearly-K\"{a}hler structure on $M$ is locally homogeneous under the action of the $7$-dimensional group $\text{Aut}_{\text{O}(2)}$.  In fact, $M$ is locally isomorphic to $\text{Aut}_{\text{O}(2)}/\text{O}(2)$, and $\text{Aut}_{\text{O}(2)}$ is a finite quotient of $\text{SU}(2) \times \text{SU}(2) \times \text{U}(1)$.  The quotient surface $\Sigma$ is locally isometric to the round sphere $\mathbb{S}^2(\frac{1}{2\sqrt{3}})$. \\
\indent (c) In cases (ii) and (iii), the nearly-K\"{a}hler structure on $M$ is cohomogeneity-one under the action of the $5$-dimensional group $\text{Aut}_{\text{O}(2)}$.  The underlying metric on $M$ is incomplete. \\
\indent In case (ii), the surface $\Sigma$ admits a Killing vector field and its Gauss curvature $K \geq 3$ lies in the interval $[ 6( c^2 - 1 - c\sqrt{c^2 - 3} ) , 6( c^2 - 1 + c \sqrt{c^2 - 3}) ]$.  In case (iii), the surface $\Sigma$ is locally isometric to the round sphere $\mathbb{S}^2(\frac{1}{\sqrt{3}})$. \\

\noindent \textit{Proof:} (a) We observe that
$$d\left( \frac{a_0}{ a_0^2 + a_1^2 + a_2^2 + 3} \right) = 0,$$
so that on a connected set
\begin{equation*}
\frac{a_0}{ a_0^2 + a_1^2 + a_2^2 + 3} = \frac{1}{2c} \tag{6.6}
\end{equation*}
where $c \in (\mathbb{R} - 0) \cup \{\infty\}$ is a constant.  Thus, the image of $A \colon P \to \mathbb{R}^3$ is a subset of: \\
\indent (i) $|c| = \sqrt{3}$: One of the points $(0, 0, \pm \sqrt{3})$. \\
\indent (ii) $\sqrt{3} < |c| < \infty$: The $2$-sphere $\{a_1^2 + a_2^2 + (a_0 - c)^2 = c^2 - 3\}$. \\
\indent (iii) $|c| = \infty$: The $2$-plane $\{a_0 = 0\}$. \\

\indent It remains to check the openness claims in (ii) and (iii).  For this, note first that either $\text{rank}(A) = 0$ or $\text{rank}(A) = 2$.  Indeed, (6.6) shows that $\text{rank}(A) \leq 2$, while an examination of the $2 \times 2$ minors in (6.2) shows that $\text{rank}(A) = 1$ is impossible. Formula (6.2) also shows that $\text{rank}(A) = 0$ if and only if the image of $A$ is one of the points $(0, 0, \pm \sqrt{3})$.  Consequently, $\text{rank}(A) = 2$ if and only if the image of $A$ is a subset of the $2$-sphere or $2$-plane described above. \\

\indent (b) In case (i), we have $(a_1, a_2, a_0) = (0,0, \pm \sqrt{3})$.  Set $(\beta_1, \beta_2) := (2\sqrt{3}\,\omega^1, \,2\sqrt{3}\,\omega^2)$ and $\beta_6 := 6\omega^6 + \sqrt{3}\,\phi$.  Then the structure equations (6.5) are:
\begin{align*}
d\beta_1 & = \beta_2 \wedge \phi & d\beta_3 & = -\beta_4 \wedge \beta_5 & d\beta_6 & = 0 \\
d\beta_2 & = \phi \wedge \beta_1 & d\beta_4 & = -\beta_5 \wedge \beta_3 \\
d\phi & = \beta_1 \wedge \beta_2 & d\beta_5 & = -\beta_3 \wedge \beta_4.
\end{align*}
Thus, $P$ with the coframing $(\omega^1, \ldots, \omega^6, \phi)$ is locally isomorphic to $\text{SU}(2) \times \text{SU}(2) \times \text{U}(1)$ with its left-invariant coframing.  In particular, $M^6$ is (locally) homogeneous.  The surface $\Sigma$ has Gauss curvature $K = 3(a_0^2 +1) = 12$, hence is locally isometric to $\mathbb{S}^2(\frac{1}{2\sqrt{3}})$. \\

\indent (c) We now consider cases (ii) and (iii).  In these cases, $\text{rank}(A) = 2$, so $H := \text{Aut}_{\text{O}(2)}$ is a $5$-dimensional group.  Since $\dim(H) = 5$, the $H$-action on $M^6$ must have cohomogeneity $\geq 1$.  On the other hand, since $H \supset G$ and $G$ acts with cohomogeneity two, the $H$-action must have cohomogeneity $\leq 2$.  If $H$ acted on $M$ with cohomogeneity two, then Proposition 4.5 would imply that $\dim(H) = 4$, a contradiction.  Thus, $M$ is cohomogeneity one under the $H$-action. \\
\indent From Podest\'{a} and Spiro's study \cite{MR2587385} of cohomogeneity-one nearly-K\"{a}hler metrics, we deduce that the underlying metric is incomplete. \\

\indent Let $\varpi \colon F^3 \to \Sigma^2$ denote the orthonormal coframe bundle of $\Sigma$, as in $\S$4.6.  In case (ii), a direct calculation using (6.1) shows that a vector field $X \in \Gamma(TF)$ satisfies $\mathcal{L}_X \omega^1 = \mathcal{L}_X \omega^2 = \mathcal{L}_X\phi = 0$ and $\mathcal{L}_X A = 0$ if and only if
$$\begin{bmatrix} \omega^1(X) \\ \omega^2(X) \\ \phi(X) \end{bmatrix} = \frac{C}{\sqrt{|a_0|}} \begin{bmatrix} -a_2 \\ a_1 \\ 3 - a_0^2 \end{bmatrix}$$
for some constant $C$.  Such a vector field $X \in \Gamma(TF)$ is $\varpi$-related to a (non-zero) Killing vector field on $\Sigma$.  Moreover, the Gauss curvature of $\Sigma$ is $K = 3(a_0^2 + 1)$, where $(a_1, a_2, a_0)$ lies on the $2$-sphere $\{a_1^2 + a_2^2 + (a_0 - c)^2 = c^2 - 3\}$, whence $a_0 \in [c - \sqrt{c^2 - 3}, c + \sqrt{c^2 - 3}]$. \\
\indent In case (iii), the surface $\Sigma$ has Gauss curvature $K = 3$, so is locally isometric to $\mathbb{S}^2(\frac{1}{\sqrt{3}})$. $\lozenge$

\subsection{An Explicit Example}

\indent \indent In case (iii), where $a_0 = 0$, we can go further and ``integrate" the structure equations (6.5) and (6.2).  Since $(a_1, a_2) \colon P \to \mathbb{R}^2$ is an $\text{O}(2)$-equivariant function, we may adapt frames as follows: consider the subset
$$P_1 := \{p \in P \colon a_2(p) = 0\}.$$
Over the open dense subset $M^* = \{m \in M \colon ((a_1)^2 + (a_2)^2)(m) \neq 0\} \subset M$, the projection $P_1 \to M^*$ is a $\mathbb{Z}_2$-bundle.  For the remainder of this section, we work on $P_1$.  For ease of notation, set $a := a_1$. \\
\indent On $P_1$, the $1$-form $\phi$ is no longer a connection form, but rather $\phi = k_1 \omega^1 + k_2 \omega^2$ for some new $G$-invariant functions $k_1, k_2 \colon P_1 \to \mathbb{R}$.  Thus, equations (6.2) now read:
\begin{align*}
da & = -(a^2 + 3)\,\omega^1 \\
0 & = -3\omega^2 + a\phi = ak_1\,\omega^1 + (ak_2 - 3)\,\omega^2,
\end{align*}
whence $(k_1, k_2) = (0, 3/a)$, and thus $\phi = \frac{3}{a}\omega^2$.  The structure equations (6.5) now read:
\begin{align*}
d\omega^1 & = 0 & d\beta_3 & = -\beta_4 \wedge \beta_5 \\
d\omega^2 & = \textstyle -\frac{3}{a}\,\omega^1 \wedge \omega^2 & d\beta_4 & = -\beta_5 \wedge \beta_3 \\
d\omega^6 & = a\,\omega^1 \wedge \omega^6 & d\beta_5 & = -\beta_3 \wedge \beta_4.
\end{align*}
Since
$$\omega^1 = -\frac{1}{a^2 + 3}\,da,$$
we see that
\begin{align*}
d\omega^2 & = d\!\left( \log \frac{a}{\sqrt{a^2 + 3}} \right) \wedge \omega^2 & d\omega^6 & = d\!\left( \log \frac{1}{\sqrt{a^2 + 3}} \right) \wedge \omega^6,
\end{align*}
and thus
\begin{align*}
\omega^2 & = \frac{a}{\sqrt{a^2 + 3}}\,dt & \omega^6 & = \frac{1}{\sqrt{a^2 + 3}}\,ds
\end{align*}
for some locally-defined functions $t,s \colon P_1 \to \mathbb{R}$.  Unwinding the notational changes, our substitution (6.3) is simply:
\begin{align*}
\beta_3 & = 3\omega^3 - \frac{3}{\sqrt{a^2 + 3}}\,dt & \omega^3 & = \frac{1}{3}\left( \beta_3 + \frac{3}{\sqrt{a^2 + 3}}\,dt \right) \\
\beta_4 & = 3\omega^4 & \omega^4 & = \frac{1}{3}\,\beta_4 \\
\beta_5 & = 3\omega^5 + \frac{a}{\sqrt{a^2 + 3}}\,ds & \omega^5 & = \frac{1}{3}\left( \beta_5 - \frac{a}{\sqrt{a^2 + 3}}\,ds \right)
\end{align*} \\
We draw the following conclusion: \\

\noindent \textbf{Theorem 6.4:} Let $(a,t,s) \colon \mathbb{R}^3 \times \mathbb{S}^3 \to \mathbb{R}^3$ denote the projection to the $\mathbb{R}^3$-factor, and let $(\beta_3, \beta_4, \beta_5)$ denote the left-invariant coframing on $\text{SU}(2) \simeq \mathbb{S}^3$.  Then the Riemannian metric $g_\infty$ on $\mathbb{R}^3 \times \mathbb{S}^3$ given by
\begin{align*}
g_\infty & = (\omega^1)^2 + (\omega^2)^2 + (\omega^6)^2 + (\omega^3)^2 + (\omega^4)^2 + (\omega^5)^2 \\
& = \frac{1}{(a^2 + 3)^2} \,da^2 + \frac{a^2}{a^2+3} dt^2 + \frac{1}{a^2+3} ds^2 \\
& \ \ \ \ + \frac{1}{9}\left(  \beta_3 + \frac{3}{\sqrt{a^2 + 3}}\,dt \right)^2 + \frac{1}{9}\beta_4^2 + \frac{1}{9}\left( \beta_5 - \frac{a}{\sqrt{a^2 + 3}}\,ds \right)^2
\end{align*}
is a Type III nearly-K\"{a}hler metric with $\mathfrak{g} = \mathfrak{su}(2) \oplus \mathfrak{u}(1)$, $u = 0$, and $a_0 = 0$.  In fact, it is the unique such metric up to diffeomorphism.

\section{Appendix}

\indent \indent In Lemma 5.7(b), a calculation using M\textsc{aple} shows that the polynomial functions $T_{41}, T_{42}$, $H_{11}, H_{12}, H_{41}, H_{41}$, and $G_1, G_2$ are:
\begin{align*}
T_{41} & = -t_1\!\left( 64t_2^3(2h_1 + 9) + 16t_2^2\ell_2 - 4t_2(h_1 -2) \right) + t_2\!\left( 6 \ell_1 - 4h_4 t_3 \right) + \ell_2 t_4 \\
T_{42} & = \textstyle 8t_2^2\left( 2t_1(\ell_1 + 12t_4) + 4t_1^2 + 6h_1 + 27 \right) - t_4(4t_1 + \ell_1 + 12 t_4)  - \frac{1}{2}\!\left(t_1^2 + t_3^2 + 3\right)  + 6\ell_2 t_2
\end{align*}
and
\begin{align*}
H_{11} & = 8t_2(t_1 \ell_1 - 4 h_1^2 - 12h_1) + 2\ell_2(h_1 + 3) \\
H_{12} & = 2t_1( 16t_2^2(2h_1 + 9) + 4 \ell_2 t_2  + h_1 + 2) - 2\ell_1( h_1 + 3) - 2t_3h_4  \\
H_{21} & = -8t_2(t_3\ell_1 + 4 h_1 h_4 + 6 h_4) + 2\ell_2 h_4 \\
H_{22} & = -2t_3( 16t_2^2(2h_1 + 9) + 4\ell_2t_2 + h_1 + 2) + 2h_4(64t_1 t_2^2 - 8t_4 - t_1 - \ell_1)
\end{align*}
and
\begin{align*}
G_1 & = 256 t_1 t_2^3 ( 2h_1 + 9)  + 8t_2(2 h_1 t_1 - 2 h_4 t_3 - 15 \ell_1 + 4 t_1) + 64 \ell_2 t_1 t_2^2 \\
G_2 & = 96t_2^2(2h_1 - 9) + 4t_2\ell_2(2 h_1 - 21) - 3(\ell_1^2 + \ell_2^2) - 4(h_1^2 + h_4^2)  - 4\ell_1 t_4 - 12h_1 - 12.
\end{align*}

\bibliographystyle{plain}
\bibliography{NKRef}

\Addresses

\pagebreak

\end{document}